\documentclass[12pt,a4paper]{amsart}
\usepackage{geometry}
\usepackage{tabls}
\usepackage{url}
\usepackage[utf8]{inputenc}
\usepackage{amsfonts}
\usepackage{amssymb}
\usepackage{mathrsfs}
\usepackage{amsmath}
\usepackage{amsthm}
\usepackage{amscd}
\usepackage{fancyhdr}
\usepackage{enumerate}
\usepackage{paralist}
\usepackage{xypic}
\usepackage{graphicx}
\usepackage{cite}
\usepackage{lipsum}
\usepackage{rotating}
\usepackage{footmisc}
\usepackage[all,cmtip]{xy}

\numberwithin{equation}{section}

\theoremstyle{plain}

\newtheorem{Def}[equation]{Definition}
\newtheorem{Thm}[equation]{Theorem}
\newtheorem{lem}[equation]{Lemma}
\newtheorem{prop}[equation]{Proposition}
\newtheorem{rem}[equation]{Remark}
\newtheorem{ex}[equation]{Example}
\newtheorem{conj}[equation]{Conjecture}

\begin{document}

\title{Unit cyclotomic multiple zeta values  for $\mu_2,\mu_3$ and $\mu_4$}

\author{Jiangtao Li}
\email{lijiangtao@csu.edu.cn}
\address{Jiangtao Li \\School of Mathematics and Statistics, HNP-LAMA,
          Central South University,
         Hunan Province, China}

\begin{abstract} Denote by $\epsilon$ a primitive root of $N^{th}$-unity. In this paper, we show that   
the unit cyclotomic multiple zeta values for $\mu_N$ generate all the cyclotomic multiple zeta values for $\mu_N$ in cases $N=2,3,4$. 
 Moreover, the unit cyclotomic multiple zeta values for $\mu_N$ can be written as $\mathbb{Q}$-linear combinations of $\left(\zeta\binom{1}{\epsilon}\right)^n, \left(\zeta\binom{1}{\epsilon^{-1}}\right)^n$ and lower depth terms in each weight $n$ in case of $N=2,3$ and $4$. By detailed analysis of the motivic Galois action, we  compute the coefficients of  $\left(\zeta\binom{1}{\epsilon}\right)^n, \left(\zeta\binom{1}{\epsilon^{-1}}\right)^n$ in the above expressions of unit cyclotomic multiple zeta values.
  \end{abstract}
\maketitle

\let\thefootnote\relax\footnotetext{
2020 $\mathnormal{Mathematics} \;\mathnormal{Subject}\;\mathnormal{Classification}$. 11M32.\\
$\mathnormal{Keywords:}$  Multiple zeta values, Cyclotomic field. \\
}

\section{Introduction}\label{int}
     For $N\geq 1$, denote by $\mu_N$ the roots of $N^{th}$-unity and  $\epsilon$ a primitive root of $N^{th}$-unity. The cyclotomic multiple zeta values for $\mu_N$ are defined by the following series:
\[
\zeta\binom{k_1,k_2,\cdots,k_r}{\epsilon_1,\epsilon_2,\cdots,\epsilon_r}=\sum_{0<n_1<n_2<\cdots<n_r}\frac{\epsilon_1^{n_1}\epsilon_2^{n_2}\cdots \epsilon_r^{n_r}}{n_1^{k_1}n_2^{k_2}\cdots n_r^{k_r}},k_i\geq 1, \epsilon_i\in \mu_N, (k_r,\epsilon_r)\neq (1,1).
\]
The condition $(k_r,\epsilon_r)\neq (1,1)$ ensures the convergence of the above series.
For cyclotomic multiple zeta value $\zeta\binom{k_1,k_2,\cdots,k_r}{\epsilon_1,\epsilon_2,\cdots,\epsilon_r}$, $K=k_1+k_2+\cdots+k_r$ is called its weight and $r$ is called its depth.
For $N=1$, they are classical multiple zeta values.

Define $\mathcal{Z}_0=\mathbb{Q}$ and $\mathcal{Z}_K$ the $\mathbb{Q}$-linear combinations of weight $K$ cyclotomic multiple zeta values for $\mu_N$. Define 
\[
\mathcal{Z}=\bigoplus_{K\geq 0} \mathcal{Z}_K.
\]
From iterated integral representations of cyclotomic multiple zeta values, it is easy to show that $\mathcal{Z}$ is a graded commutative $\mathbb{Q}$-algebra.

Cyclotomic multiple zeta values have been studied by Deligne, Goncharov, Hoffman, Racinet, Zhao,$\cdots$  in a series of papers.

Brown \cite{brown} introduced the definition of motivic multiple zeta values. By detailed analysis of the motivic Galois action on motivic multiple zeta values.  Brown proved  a conjecture of Hoffman \cite{hoff}.

Glanois \cite{gla} introduced the definitions of cyclotomic motivic multiple zeta values for $\mu_N$, $N=2,3,4,6,8$. 
Glanois gave a basis of cyclotomic motivic multiple zeta values in each case respectively.

Denote by $\mathcal{Z}^{(1)}_K$ the $\mathbb{Q}$-linear space generated by the following
weight $K$ elements:
\[
\zeta\binom{1,\,1,\,\cdots,\,1}{\epsilon_1,\epsilon_2,\cdots,\epsilon_K},\epsilon_i\in \mu_N, \epsilon_K\neq 1.
\]
We call  $\mathcal{Z}^{(1)}_K$ unit cyclotomic multiple zeta values of weight $K$.
 Define 
\[
\mathcal{Z}^{(1)}=\bigoplus_{K\geq 0} \mathcal{Z}^{(1)}_K.
\]
Clearly $\mathcal{Z}^{(1)}$ is a  graded $\mathbb{Q}$-subalgebra of $\mathcal{Z}$.
 Unit cyclotomic multiple zeta values have been studied by  Borwein, Bradley, Broadhurst and Lisonek \cite{BBBL} for $N=2$. They have also been studied by Zhao \cite{zhao} for $N=3, 4$. Zhao  proposed the following conjecture (Conjecture $6.2$ in \cite{zhao}).
 \begin{conj}\label{zhao}
 The set of following elements 
\[
\zeta\binom{1,\;1,\;\cdots,\;1}{\epsilon_1,\epsilon_2,\cdots,\epsilon_K}, \epsilon_i\in \{\epsilon, \epsilon^2\}
\]
is a basis for cyclotomic multiple zeta values of weight $K$ for $N=3$ and $4$ for  $K\geq 1$. 
\end{conj}

In this paper, it is proved  that
\begin{Thm}\label{uzhao}
For $N=2,3,4$, all the cyclotomic multiple zeta values for $\mu_N$ are $\mathbb{Q}$-linear combinations of unit cyclotomic multiple zeta values for $\mu_N$.
\end{Thm}
Theorem \ref{uzhao} provides partial evidence towards Zhao's conjecture. For a multiple zeta value $\zeta(k_1,\cdots,k_r)$,   Ohno and Zagier \cite{OZ} defined the height  of $\zeta(k_1,\cdots,k_r)$ as $$i=\#\big{\{} k_j\,\big{|}\,k_j>1\big{\}}.$$
Brown's main result \cite{brown} implies that multiple zeta values are generated by multiple zeta values of {\bf maximal height}.
In the cyclotomic cases, Theorem \ref{uzhao} shows that the cyclotomic multiple zeta values for $\mu_2,\mu_3,\mu_4$ are generated by cyclotomic multiple zeta values of {\bf minimal height} (which are unit cyclotomic multiple zeta values.)

In Section \ref{mainsection} we will give a set of generating elements  (which are $\mathbb{Q}$-linear combinations of unit cyclotomic multiple zeta values for $\mu_N$) of cyclotomic multiple zeta values for $\mu_N$ in cases $N=2,3,4$.

Denote by $\mathcal{D}_{r}\mathcal{Z}_K$ the $\mathbb{Q}$-linear space generated by weight $K$ and depth $\leq r$ cyclotomic multiple zeta values for $\mu_N$. Define 
\[
gr_r^{\mathcal{D}}\mathcal{Z}_K=\mathcal{D}_r\mathcal{Z}_K/\mathcal{D}_{r-1}\mathcal{Z}_K.
\]
Denote by $gr_r^{\mathcal{D}}\mathcal{Z}^{(1)}_r$ the $\mathbb{Q}$-linear subspace of $gr_r^{\mathcal{D}}\mathcal{Z}_r$ which is generated by the images of unit cyclotomic multiple zeta values of weight $r$ and depth $r$.

In this paper, by generalizing the motivic method of Brown \cite{depth} to the cyclotomic case, we will study the structure of $gr_r^{\mathcal{D}}\mathcal{Z}^{(1)}_r$ in each weight $r$ for  $\mu_2, \mu_3$ and $\mu_4$.

\begin{Thm}\label{gene}
(i) For $N=2$, $gr_r^{\mathcal{D}}\mathcal{Z}_r^{(1)}$ is generated by the image of 
\[
\zeta\dbinom{1,\cdots,1,\;\;1\,}{\underbrace{1,\cdots,1,-1}_r}
\]
in $gr_r^{\mathcal{D}}\mathcal{Z}_K$ as a $\mathbb{Q}$-linear subspace.\\
(ii) For $N=3,4$, $gr_r^{\mathcal{D}}\mathcal{Z}_r^{(1)}$ is generated by the images of
\[
\zeta\dbinom{1,\cdots,1,1}{\underbrace{1,\cdots,1,\epsilon}_r}, \zeta\dbinom{1,\cdots,1,\,\;1\;}{\underbrace{1,\cdots,1,\epsilon^{-1}}_r}
\]
in $gr_r^{\mathcal{D}}\mathcal{Z}_r$ as a $\mathbb{Q}$-linear subspace.
\end{Thm}

The essential reason behind Theorem \ref{gene} is that most parts of the motivic Galois action on the motivic version of $gr_r^{\mathcal{D}}\mathcal{Z}_K^{(1)}$ vanish. As a result, the motivic version of $gr_r^{\mathcal{D}}\mathcal{Z}_K^{(1)}$ is just a linear subspace of dimension one or two.

From the iterated integral representation of cyclotomic multiple zeta values, it is easy to check that 
\[
\zeta\dbinom{1,\cdots,1,1}{\underbrace{1,\cdots,1,\epsilon}_r}=\frac{1}{r!}\left[\zeta\binom{1}{\epsilon} \right]^r=\frac{(-1)^r}{r!}\left[\mathrm{log}\;(1-\epsilon) \right]^r.
\]
Thus for 
\[
\zeta\binom{1,\;1,\;\cdots,1}{\epsilon_1,\epsilon_2,\cdots, \epsilon_r}, \epsilon_i\in \mu_N, N=2,3,4,
\]
we have 
\[
\zeta\binom{1,\;1,\;\cdots,1}{\epsilon_1,\epsilon_2,\cdots, \epsilon_r}=c_{\epsilon_1,\cdots,\epsilon_r}\left(\mathrm{log}\;2\right)^r+\mathrm{LDT}, \;\epsilon_i\in \{\pm1\}, \epsilon_r=-1
\]
and for $N=3,4$,
\[
\zeta\binom{1,\;1,\;\cdots,1}{\epsilon_1,\epsilon_2,\cdots, \epsilon_r}=a_{\epsilon_1,\cdots,\epsilon_r}\left[\mathrm{log}\;(1-\epsilon)\right]^r+b_{\epsilon_1,\cdots,\epsilon_r}\left[\mathrm{log}\;(1-\epsilon^{-1})\right]^r+\mathrm{LDT},
\]
where $a_{\epsilon_1,\cdots,\epsilon_r}, b_{\epsilon_1,\cdots,\epsilon_r}, c_{\epsilon_1,\cdots,\epsilon_r}\in \mathbb{Q}$. Here the symbol LDT is short for lower depth terms.

In fact, for the numbers $$a_{\epsilon_1,\cdots,\epsilon_r}, \; b_{\epsilon_1,\cdots,\epsilon_r}, \;c_{\epsilon_1,\cdots,\epsilon_r},$$ we have
\begin{Thm}\label{efca}
$(i)$ For $N=2$, $\eta_1,\cdots,\eta_{r-1}\in \{\pm 1\}, \eta_r=-1$, one has
\[
\zeta\binom{1,\;1,\;\cdots,1}{\eta_1,\eta_2,\cdots, \eta_r}=\zeta\binom{1,\cdots,1,\;\;1}{\underbrace{1,\cdots,1,-1}_r}+\mathrm{LDT}.\]
$(ii)$ For $N=3$, $\epsilon_1,\cdots,\epsilon_{r-1}\in \mu_3, \epsilon_r\neq 1$, one has
\[
\zeta\binom{1,\;1,\;\cdots,1}{\epsilon_1,\epsilon_2,\cdots, \epsilon_r}=\zeta\binom{1,\cdots,1,\,1}{\underbrace{1,\cdots,1,\epsilon_r}_r}+\mathrm{LDT}.\]
$(iii)$ For $N=4$, one has
\[
\zeta\binom{1,\cdots,1,\;\;1}{\underbrace{1,\cdots,1,-1}_r }=2^{r-1}\left(   \zeta\binom{1,\cdots,1,1}{\underbrace{1,\cdots,1,\epsilon}_r }+    \zeta\binom{1,\cdots,1,\;1}{\underbrace{1,\cdots,1,\epsilon^{-1}}_r }     \right)+\mathrm{LDT}.\tag{1}\]
\[\zeta\binom{1,\;1,\;\cdots,1}{\eta_1,\eta_2,\cdots, \eta_s}=\zeta\binom{1,\cdots,1,\,1\;}{1,\cdots,1,\eta_s}+\mathrm{LDT},\;  \eta_j \in \{\pm 1\}, \forall\, j ,\eta_s=-1,\tag{2}
\]
\[
\zeta\binom{1,\;1,\;\cdots,1\,}{\epsilon_1,\epsilon_2,\cdots, \epsilon_r}=
\zeta\binom{1,\cdots,1,\,1\,}{1,\cdots,1,\epsilon_r}+\mathrm{LDT},\; \epsilon_j \in \{\epsilon,\epsilon^{-1}\}, \forall\,j.\tag{3}
\]
Furthermore, for $\epsilon_1,\cdots, \epsilon_r\in \{\epsilon,\epsilon^{-1}\}, \eta_1,\cdots, \eta_s\in \{\pm 1\}, \eta_1=-1, r,s\geq 1$, one has
\[\small
\begin{split}
&\;\;\;\;\zeta\binom{1,\,\cdots,1,\;1\,,\,\cdots,1}{\epsilon_r,\cdots,\epsilon_1,\eta_s,\cdots,\eta_1}\\
&=\sum_{j=0}^{s-1}\binom{r+s-1}{j}\cdot \zeta\binom{1,\cdots,1,\,1\,\;}{\underbrace{1,\cdots,1,\eta_1}_{r+s}}+\binom{r+s-1}{s}\zeta\binom{1,\cdots,1,\,\;\;\;\;\;\;\,1\;\;\;\;\;\;}{\underbrace{1,\cdots,1,\epsilon_1\eta_1\cdots \eta_s}_{r+s}}+\mathrm{LDT}.
\\
\end{split}\tag{4}
\]
For $\epsilon_1\in \{\epsilon,\epsilon^{-1}\}, \eta_1,\cdots, \eta_s\in \{\pm 1\}, s\geq 1$, one has
\[
\begin{split}
&\;\;\;\;\zeta\binom{1\,,\,\cdots,1,\;1\,}{\eta_s,\cdots,\eta_1,\epsilon_1}\\
&=2^s\zeta\binom{ 1,\cdots,1,\,1\,}{\underbrace{1,\cdots,1}_s, \epsilon_1}+\sum_{j=1}^s(-1)^j2^{s-j}\binom{s}{j}\zeta\binom{1,\cdots,1,\,\;\;\;\;\;\;\,1\,\;\;\;\;\;\;}{\underbrace{1,\cdots,1}_s,\underbrace{\eta_j\cdots \eta_1}_j \epsilon_1}+\mathrm{LDT}. \\ 
\end{split} \tag{5}
\]
For $\epsilon_1,\cdots , \epsilon_r\in \{\epsilon,\epsilon^{-1}\}, \eta_1,\cdots, \eta_s\in \{\pm 1\}, r\geq 2, s\geq 1$, one has
\[
\begin{split}
&\;\;\;\;\zeta\binom{1\,,\,\cdots,1\,,\;1\,,\cdots, \,1}{\eta_s,\cdots,\eta_1,\epsilon_r,\cdots,\epsilon_1}=\zeta\binom{1,\cdots,1,\,1\;}{\underbrace{1,\cdots,1,\epsilon_1}_{r+s}}+\mathrm{LDT}.\\
\end{split} \tag{6}
\]
\end{Thm}

\section{Mixed Tate motives}
In this section we will give a brief introduction to mixed Tate motives. For more details, see \cite{GF}, \cite{del}  and \cite{deligne}. 
Since we only discuss cyclotomic multiple zeta values for $N=2,3$ and $4$, the number $N$ in this section is $2,3$ or $4$.

\subsection{Mixed Tate motives over $\mathcal{O}_N[\frac{1}{N}]$}\label{mtm}

     Denote by $\mathcal{O}_N$ the algebraic integer ring of the cyclotomic field $\mathbb{Q}[\mu_N]$. Deligne and Goncharov \cite{deligne} constructed the category of mixed Tate motives over $\mathcal{O}_N[\frac{1}{N}]$. Denote it by $\mathcal{MT}_N$ for short. $\mathcal{MT}_N$ is a neutral Tannakian catogory with the natural fiber functor
\[
\omega:\mathcal{MT}_N\rightarrow \mathrm{Vect}_{\mathbb{Q}}; M\mapsto \bigoplus \omega_r(M),
\]
where
\[
\omega_r(M)=\mathrm{Hom}_{\mathcal{MT}_N}(\mathbb{Q}(r),gr_{-2r}^{\omega}(M)).
\]

Let $\mathcal{G}^{\mathcal{MT}_2}$ be the Tannakian fundamental group of $\mathcal{MT}_N$  under this fiber functor, then we have 
$$\mathcal{G}^{\mathcal{MT}_{N}}=\mathbb{G}_{m}\ltimes \mathcal{U}^{\mathcal{MT}_{N}},$$
where $\mathcal{U}^{\mathcal{MT}_N}$ is a pro-unipotent algebraic group.

From Deligne and Goncharov's construction \cite{deligne} and Borel's theorem on K-group  of number fields, we have 
\[
\mathrm{Ext}^1_{\mathcal{MT}_N}(\mathbb{Q}(0),\mathbb{Q}(n))\cong \mathbb{Q}, \mathrm{if}\;N=2,3,4,n\geq 1,\mathrm{odd},
\]
\[
\mathrm{Ext}^1_{\mathcal{MT}_N}(\mathbb{Q}(0),\mathbb{Q}(n))= 
\begin{cases}
0, &\mathrm{if}\; N=2,n\leq 0 \;\mathrm{or}\; n\; \mathrm{even},\\
\mathbb{Q},&\mathrm{if}\; N=3,4, n\geq 2,\mathrm{even},\\
0,&\mathrm{if}\; N=3,4, n\leq 0,\\
\end{cases}
\]
\[
\mathrm{Ext}^2_{\mathcal{MT}_N}(\mathbb{Q}(0),\mathbb{Q}(n))=0, \forall\,n \in \mathbb{Z}.
\]

Denote by $\mathfrak{g}_N$ the Lie algebra of $\mathcal{U}^{\mathcal{MT}_N}$. From the above facts about extension groups, we know that $\mathfrak{g}_N$ is a free Lie algebra. Its generators are $\sigma_{2n+1},n\geq 0$ (weight $\sigma_{2n+1}=-2n-1$) for $N=2$ and $\sigma_{n},n\geq 1$ for $N=3,4$.

From the natural correspondence between pro-nilpotent Lie algebra and pro-unipotent group, we have that 
\[
\mathcal{O}(\mathcal{U}^{MT_N})\cong\begin{cases}
\mathbb{Q}\langle f_1, f_3,\cdots, f_{2n+1}\cdots\rangle,&N=2,\\
\mathbb{Q}\langle f_1,f_2,\cdots, f_n,\cdots\rangle, &N=3,4,\\
\end{cases}
\]
as a graded $\mathbb{Q}$-algebra, where the multiplication on the right side is actually the shuffle product $\rotatebox{90}{$\rotatebox{180}{$\exists$}$}$ on the non-commutative  word sequences in  $f_n, n\geq 1$. It is given by the following induction formulas:
\[
1\,\rotatebox{90}{$\rotatebox{180}{$\exists$}$}\,w=w\,\rotatebox{90}{$\rotatebox{180}{$\exists$}$}\,1=w,
\]
\[
uw_1\, \rotatebox{90}{$\rotatebox{180}{$\exists$}$}\, vw_2=u(w_1\,\rotatebox{90}{$\rotatebox{180}{$\exists$}$}\,vw_2)+v(uw_1\,\rotatebox{90}{$\rotatebox{180}{$\exists$}$}\,w_2),
\]
where $u,v \in \{f_n, n\geq 1\}$. In fact, $f_n,n\geq 1$ are dual to $\sigma_n,n\geq 1$ in the natural way. 

\subsection{Motivic cyclotomic multiple zeta values}

From \cite{deligne}, the motivic fundamental groupoid of $\mathbb{P}^1-\{0,\mu_N,\infty\}$ can be realized in the category $\mathcal{MT}\left(\mathcal{O_N}[\frac{1}{N}]\right)$. Denote by ${}_0\Pi_1$ the motivic fundamental groupoid of $\mathbb{P}^1-\{0,\mu_N,\infty\}$ from $\overrightarrow{1}_0$ to $\overrightarrow{-1}_1$ (the tangential vector $\overrightarrow{1}$ at point $0$ and the tangential vector $\overrightarrow{-1}$ at the point $1$). Its  ring of regular functions is isomorphic to
\[
\mathcal{O}({}_0\Pi_1)\cong (\mathbb{Q}\langle e^0, e^{\mu_N}\rangle,\rotatebox{90}{$\rotatebox{180}{$\exists$}$} )
\]
under the Tannakian correspondence, where $\mathbb{Q}\langle e^0,e^{\mu_N}\rangle$ is the non-commutative polynomial linear space in the words $e^0,e^{\epsilon},\epsilon\in \mu_N$ with the shuffle product $\rotatebox{90}{$\rotatebox{180}{$\exists$}$}$ (its definition is similar to the one in Section \ref{mtm}). As a result, $(\mathbb{Q}\langle e^0, e^{\mu_N}\rangle,\rotatebox{90}{$\rotatebox{180}{$\exists$}$})$ is a commutative $\mathbb{Q}$-algebra. Under Tannakian correspondence, the ring of regular functions of $\mathcal{U}^{\mathcal{MT}_N}$ has a coaction on $\mathcal{O}({}_0\Pi_1)$.

For arbitrary word sequence $u_1u_2\cdots u_k$ in $e^0,e^{\mu_N}$, if $\delta,\eta \rightarrow 0$, by direct calculation it is easy to check that (see the Appendix A in 
\cite{LM})
\[
\mathop{\int\cdots\int}\limits_{\delta<t_1<\cdots<t_k<1-\eta}\omega_{u_1}(t_1)\cdots \omega_{u_k}(t_k)=P(\mathrm{log}(\delta),\mathrm{log}(\eta))+O\left(\mathrm{sup}(\delta|\mathrm{log}(\delta)|^A+\eta|\mathrm{log}(\eta)|^B)\right),
\]
where $\omega_{e^0}(t)=\frac{dt}{t},\omega_{e^{\epsilon}}(t)=\frac{dt}{\epsilon-t}$ for $\epsilon\in \mu_N$ and $P$ is a $\mathbb{C}$-coefficients polynomial of two variables.

Define $dch:\mathcal{O}({}_0\Pi_1)=\mathbb{Q}\langle e^0,e^{\mu_N}\rangle\rightarrow \mathbb{C}$ by 
\[
dch(u_1u_2\cdots u_k)=P(0,0).
\]
One can check that the images of $\mathcal{O}({}_0\Pi_1)$
are $\mathbb{Q}$-linear combinations of cyclotomic multiple zeta values (see also the Appendix A in 
\cite{LM}). By the shuffle product of the iterated integrals,  $dch$ is a ring homomorphism 
\[
dch:\mathcal{O}({}_0\Pi_1)=\mathbb{Q}\langle e^0,e^{\mu_N}\rangle\rightarrow \mathbb{C}.
\]
So it also corresponds to a point $dch\in {}_0\Pi_1(\mathbb{C})$. This point $dch$ essentially comes from the comparison isomorphism between Betti fundamental groupoid of $\mathbb{P}^1-\{0,\mu_N,\infty\}$ and de-Rham fundamental groupoid of $\mathbb{P}^1-\{0,\mu_N,\infty\}$.

Denote by $\mathcal{I}\subseteq \mathcal{O}({}_0\Pi_1)$ the kernel of $dch$. Define $J^{\mathcal{MT}}\subseteq \mathcal{I}$ the largest graded sub-ideal of $\mathcal{I}$ which is stable under the coaction of $\mathcal{O}(\mathcal{U}^{\mathcal{MT}_N})$. The motivic cyclotomic multiple zeta algebra $\mathcal{H}$ for $\mu_N$ is $\mathcal{O}({}_0\Pi_1)/\mathcal{J}^{\mathcal{MT}_N}$.

Denote by $I^{\mathfrak{m}}$ the natural quotient map
\[
I^{\mathfrak{m}}:\mathcal{O}({}_0\Pi_1)=\mathbb{Q}\langle e^0, e^{\mu_N}\rangle\rightarrow \mathcal{H}
\]
 and $per$ the map $per:\mathcal{H}\rightarrow \mathbb{C}$ satisfying $per\circ I^{\mathfrak{m}}=dch$.

The motivic multiple zeta value $\zeta^{\mathfrak{m}}\binom{n_1,n_2,\cdots,n_r}{\epsilon_1,\,\epsilon_2,\,\cdots,\, \epsilon_r}$ is defined by
\[
\zeta^{\mathfrak{m}}\binom{n_1,n_2,\cdots,n_r}{\epsilon_1,\,\epsilon_2,\,\cdots,\, \epsilon_r}=I^{\mathfrak{m}}\left(e^{(\epsilon_1\cdots \epsilon_r)^{-1}}(e^0)^{n_1-1}e^{(\epsilon_2\cdots \epsilon_r)^{-1}}(e^0)^{n_2-1}\cdots e^{\epsilon_r^{-1}}(e^0)^{n_r-1} \right).
\]
By direct calculation of the iterated integral, we have 
\[
per\left( \zeta^{\mathfrak{m}}\binom{n_1,n_2,\cdots,n_r}{\epsilon_1,\,\epsilon_2,\,\cdots,\, \epsilon_r} \right)=\zeta \binom{n_1,n_2,\cdots,n_r}{\epsilon_1,\,\epsilon_2,\,\cdots,\, \epsilon_r}
\]
for $(n_r,\epsilon_r)\neq (1,1)$.

We will need the following lemma to study the unit cyclotomic multiple zeta values:
\begin{lem}\label{cont}
The images of the elements $e^{\epsilon_1}e^{\epsilon_2}\cdots e^{\epsilon_r},\epsilon_i\in \mu_N$ in $\mathcal{O}({}_0\Pi_1)$ under the map $dch$ are elements of $\mathcal{Z}_r^{(1)}$.
\end{lem}
 \noindent{\bf Proof}:
For word sequence $e^{\epsilon_1}e^{\epsilon_2}\cdots e^{\epsilon_r}$, if $\epsilon_r\neq 1$, then the integral
\[
\mathop{\int\cdots\int}\limits_{\delta<t_1<\cdots<t_r<1-\eta}\omega_{e^{\epsilon_1}}(t_1)\cdots \omega_{e^{\epsilon_r}}(t_r)
\]
converges when $\delta,\eta \rightarrow 0$.
So if $\epsilon_r\neq 1$, then
\[
\begin{split}
&\;\;\;\;\;dch(e^{\epsilon_1}e^{\epsilon_2}\cdots e^{\epsilon_r})\\
&= \mathop{\int\cdots\int}\limits_{0<t_1<\cdots<t_r<1}\omega_{e^{\epsilon_1}}(t_1)\cdots \omega_{e^{\epsilon_r}}(t_r) \\
&= \mathop{\int\cdots\int}\limits_{0<t_1<\cdots<t_r<1}\left(\sum_{n_1\geq 0}t_1^{n_1}\epsilon_1^{-n_1-1}\right)dt_1\cdots \left(\sum_{n_r\geq 0}t_r^{n_r}\epsilon_r^{-n_r-1}\right)dt_r\\
&=\sum_{0<n_1<n_2<\cdots<n_{r-1}<n_r}\frac{(\frac{\epsilon_2}{\epsilon_1})^{n_1}(\frac{\epsilon_3}{\epsilon_2})^{n_2}\cdots (\frac{\epsilon_r}{\epsilon_{r-1}})^{n_{r-1}}(\frac{1}{\epsilon_r})^{n_r}}{n_1n_2\cdots n_{r-1}n_r}\\
&=\zeta\dbinom{n_1,n_2,\cdots,n_{r-1},n_r}{\frac{\epsilon_2}{\epsilon_1},\frac{\epsilon_3}{\epsilon_2},\cdots,\frac{\epsilon_r}{\epsilon_{r-1}},\frac{1}{\epsilon_r}}.
\end{split}
\]

By definition we have $dch(e^1)=0$. From the shuffle product on iterated integrals, we have 
\[
\begin{split}
&\;\;\;\;\;dch(e^{\epsilon_1}e^{\epsilon_2}\cdots e^{\epsilon_r})\cdot dch(e^1)\\
&=dch(e^{\epsilon_1}e^{\epsilon_2}\cdots e^{\epsilon_r}\rotatebox{90}{$\rotatebox{180}{$\exists$}$} \;e^1) \\
&=dch(e^1e^{\epsilon_1}e^{\epsilon_2}\cdots e^{\epsilon_r}+e^{\epsilon_1}e^1e^{\epsilon_2}\cdots e^{\epsilon_r}+\cdots+ e^{\epsilon_1}e^{\epsilon_2}\cdots e^1 e^{\epsilon_r}+e^{\epsilon_1}e^{\epsilon_2}\cdots e^{\epsilon_r}e^1) \\
&=0.
\end{split}
\]
So 
\[
\begin{split}
&\;\;\;\;\;dch(e^{\epsilon_1}e^{\epsilon_2}\cdots e^{\epsilon_r}e^1)\\
&=-dch(e^1e^{\epsilon_1}e^{\epsilon_2}\cdots e^{\epsilon_r})-dch(e^{\epsilon_1}e^1e^{\epsilon_2}\cdots e^{\epsilon_r})-\cdots-dch(e^{\epsilon_1}e^{\epsilon_2}\cdots e^1 e^{\epsilon_r}).\\
\end{split}
\]
As a result, 
\[
dch( e^{\epsilon_1}e^{\epsilon_2}\cdots e^{\epsilon_r})\in \mathcal{Z}_r^{(1)},\forall \epsilon_i\in \mu_N, 1\leq i\leq r
\]
by induction.   $\hfill\Box$\\

Denote by $\mathcal{H}^{(1)}$ the images of $\mathbb{Q}\langle e^{\mu_N}\rangle$ (viewed as a $\mathbb{Q}$-subalgebra of $\mathcal{O}({}_0\Pi_1)$) under the quotient map $I^{\mathfrak{m}}:\mathcal{O}({}_0\Pi_1)=\mathbb{Q}\langle e^0, e^{\mu_N}\rangle\rightarrow \mathcal{H}$ and also denote by $\mathcal{H}^{(1)}_r$ its weight $r$ part.  By Lemma \ref{cont} we have 
\[
per(\mathcal{H}^{(1)})=\mathcal{Z}^{(1)}.
\]

In $\mathcal{O}({}_0\Pi_1)=\mathbb{Q}\langle e^0, e^{\mu_N}\rangle$, for any word $u_1\cdots u_k$, $u_i \in \{e^0,e^{\mu_N}\}$, $k$ is called its weight and the total number of occurrences of $e^{\epsilon},\epsilon\in \mu_N$ is called its depth. 
Denote by $\mathcal{D}_r\mathbb{Q}\langle e^0,e^{\mu_N}\rangle$ the subspace which consists of elements of depth $\leq r$.

From Section $6$, \cite{deligne} it follows that the depth filtration on $\mathcal{O}({}_0\Pi_1)$ is motivic. So it induces a natural depth filtration on $\mathcal{H}$. By direct calculation one can show that 
\[
per(\mathcal{D}_r\mathcal{H})=\mathcal{D}_r\mathcal{Z},\forall r\geq 0.
\]

Denote by $gr_r^{\mathcal{D}}\mathcal{H}=\mathcal{D}_r\mathcal{H}/\mathcal{D}_{r-1}\mathcal{H}$, and  define $gr_r^{\mathcal{D}}\mathcal{H}_r^{(1)}$ the natural images of weight $r$ unit cyclotomic motivic multiple zeta values  $\mathcal{H}_r^{(1)}$  in $gr_r^{\mathcal{D}}\mathcal{H}$.
In this paper we will focus on the structure of $gr_r^{\mathcal{D}}\mathcal{H}_r^{(1)}$ for all $r\geq 1$.

\subsection{Motivic Galois action}\label{mga}

In this subsection we will explain the depth-graded version motivic Galois action on the motivic cyclotomic multiple zeta values.

For $x,y\in \{0,\mu_N\}$, denote by ${}_x\Pi_y$ the motivic fundamental groupoid from  the tangential point at $x$ to the tangential point at $y$.

Under Tannakian correspondence, $\mathcal{O}({}_x\Pi_y)\cong (\mathbb{Q}\langle e^0,e^{\mu_N}\rangle,\rotatebox{90}{$\rotatebox{180}{$\exists$}$})$ for $x,y\in \{0,\mu_N\}$. There is a natural $\mu_N$-action on these groupoids: for $\epsilon\in \mu_N$, we have a morphism of schemes 
\[
\epsilon:{}_x\Pi_y\rightarrow {}_{\epsilon x}\Pi_{\epsilon y}
\]
which is defined by 
\[
\epsilon^{*}:\mathcal{O}({}_{\epsilon x}\Pi_{\epsilon y}) \rightarrow \mathcal{O}({}_x\Pi_y); e^{\alpha}\mapsto e^{\epsilon^{-1} \alpha}, \forall \alpha \in \{0,\mu_N\}
\]
on the homomorphism between rings of regular functions.

Let $V_N$ be a subgroup of automorphisms of the motivic fundamental groupoids (all base points are tangential points at $\{0,\mu_N\}$) of $\mathbb{P}^1-\{0,\mu_N,\infty\}$ satisfying the following properties:\\
(i) Elements of $V_N$  are compatible with the composition law on the motivic fundamental  groupoids of $\mathbb{P}^1-\{0,\mu_N,\infty\}$;\\
(ii) Elements of $V_N$ fix $\mathrm{exp}(e_i)\in {}_i\Pi_i$ for $i\in \{0,\mu_N\}$;\\
(iii) Elements of $V_N$ are equivariant with the $\mu_N$-action on the motivic fundamental groupoids.

By proposition 5.11 in \cite{deligne}, the following map 
\[
\xi:V_N\rightarrow {}_0\Pi_1, a\mapsto a({}_01_1)
\]
is an isomorphism of schemes and 
\[
\mathrm{Lie}\; V_N=(\mathbb{L}(e_0,e_{\mu_N}),\{\,,\,\}).
\]
Here $\mathbb{L}(e_0,e_{\mu_N})$ is the free Lie algebra generated by the symbols $e_0, e_{\epsilon},\epsilon\in \mu_N$ and $\{\;,\;\}$ denotes the Ihara Lie bracket on $\mathbb{L}(e_0,e_{\mu_N})$.

The action of $\mathcal{U}^{\mathcal{MT}_N}$ on ${}_x\Pi_y,x,y\in \{0,\mu_N\}$ factors through $V_N$. As a result, there is a Lie algebra homomorphism:
\[
i:\mathfrak{g}_N\rightarrow \mathrm{Lie}\; V_N=\left(\mathbb{L}(e_0,e_{\mu_N}),\{\;,\;\} \right).
\]
The map $i$ is injective by the main results of Deligne \cite{del} for $N=2,3,4$.

For any element $w$ in $\mathbb{L}(e_0,\mu_N)$, let $depth(w)$ be the smallest number of total occurrences of $e_{\epsilon},\epsilon\in \mu_N$ in $w$, it induces a depth decreasing filtration $\mathcal{D}$ on $\mathbb{L}(e_0,e_{\mu_N})$:
\[
\mathcal{D}^r\mathbb{L}(e_0,e_{\mu_N})=\{w\in \mathbb{L}(e_0,e_{\mu_N});depth(w)\geq r\}.
\]

We write $E^{(n)}_{\epsilon}=\mathrm{ad}(e_0)^{n}e_{\epsilon}$ for short, $\forall \epsilon\in \mu_N$.
According to Section $3.11$ in  \cite{del}, for $N=2$, the map $i$ satisfies:
\[
i(\sigma_1)=e_{-1},\eqno{(1)}
\]
\[
i(\sigma_{2n+1})=(1-2^{2n})E^{(2n)}_{-1}+2^{2n}E^{(2n)}_1+\mathrm{HDT},\forall n\geq 1.\eqno{(2)}
\]
For $N=3$, the map $i$ satisfies:
\[i(\sigma_1)=e_{\epsilon}+e_{\epsilon^{-1}},\eqno{(3)}
\]
\[
i(\sigma_{2n})=E^{(2n-1)}_{\epsilon}-E^{(2n-1)}_{\epsilon^{-1}}+\mathrm{HDT},\forall n\geq 1,\eqno{(4)}
\]
\[
i(\sigma_{2n+1})=(1-3^{2n})\left[E^{(2n)}_{\epsilon}+E^{(2n)}_{\epsilon^{-1}}\right]+2\cdot 3^{2n}E^{(2n)}_1+\mathrm{HDT},\forall n\geq 1.\eqno{(5)}
\]
For $N=4$, the map $i$  satisfies:
\[
i(\sigma_1)=e_{\epsilon}+e_{\epsilon^{-1}}+2e_{-1},\eqno{(6)}
\]
\[
i(\sigma_{2n})=E^{(2n-1)}_{\epsilon}-E^{(2n-1)}_{\epsilon^{-1}}+\mathrm{HDT},\eqno{(7)}
\]
\[
i(\sigma_{2n+1})=(1-2^{2n})\left[E^{(2n)}_{\epsilon}+E^{(2n)}_{\epsilon^{-1}}\right]+2\cdot 2^{2n}\left(1-2^{2n} \right)E^{(2n)}_{-1}+2\cdot 2^{4n}E^{(2n)}_1+\mathrm{HDT}.\eqno{(8)}
\]
In the above formulas, HDT means the higher depth terms.

The motivic Lie algebra $\mathfrak{g}_N$ has an induced depth filtration $\mathcal{D}^r\mathfrak{g}_N$ from the injective map $i$. Since Ihara bracket is compatible with the depth filtration, we know that the depth-graded space
\[
\mathfrak{dg}_N=\bigoplus_{r\geq 1}\mathcal{D}^r\mathfrak{g}_N/\mathcal{D}^{r+1}\mathfrak{g}_N
\]
is a Lie algebra with induced Ihara Bracket. By \cite{del}, $\mathfrak{dg}_N$ is a free Lie algebra for $N=2,3,4$ with generators $\overline{i(\sigma_{2n-1})}, n\geq 1$ for $N=2$ and with generators $\overline{i(\sigma_{n})},n\geq 1$ for $N=3,4$, where the symbol $\overline{i(\sigma_n})$ means the depth one parts of $i(\sigma_n)$.

The action of $\mathrm{Lie}\;V$ on $\mathcal{O}({}_0\Pi_1)$ is compatible with the depth filtration. Since the expression of $i(\sigma_{2n+1})$ in $(\mathbb{L}(e_0,e_1,e_{-1}),\{\;,\,\})$ has canonical depth one parts, for $n\geq 0$, $\sigma_{2n+1}$ in $\mathfrak{g}_2=\mathrm{Lie}\;\mathcal{U}^{\mathcal{MT}_2}$ induces a well-defined derivation
\[
\partial_{2n+1}:gr_r^{\mathcal{D}}\mathcal{H}\rightarrow gr_{r-1}^{\mathcal{D}}\mathcal{H}.
\]
For $N=3,4, n\geq 1$, $\sigma_n$ in $\mathfrak{g}_N=\mathrm{Lie}\;\mathcal{U}^{\mathcal{MT}_N}$ also induces a derivation similarly
\[
\partial_{n}:gr_r^{\mathcal{D}}\mathcal{H}\rightarrow gr_{r-1}^{\mathcal{D}}\mathcal{H}.
\]

 The explicit calculation of these derivations is very complicated.  We now give the key idea to calculate these derivations explicitly, which is essentially the generalization of Brown's observation in \cite{depth}.

Since $\mathcal{O}({}_0\Pi_1)$ is an ind-object in the category $\mathcal{MT}_N$, under Tannakian correspondence there is an action of the motivic Lie algebra
\[
\mathfrak{g}_N\times \mathcal{O}({}_0\Pi_1)\rightarrow \mathcal{O}({}_0\Pi_1).
\]

Denote by $\mathfrak{h}_N=\mathrm{Lie}\;V_N=\left(\mathbb{L}(e_0,e_{\mu_N}),\{\;,\,\} \right)$. The action of $\mathfrak{g}_N$ on $\mathcal{O}({}_0\Pi_1)$ factors through the action of $\mathfrak{h}_N$ on $\mathcal{O}({}_0\Pi_1)$.

Denote by $\mathcal{U}\mathfrak{h}_N$ the universal enveloping algebra of $\mathfrak{h}_N$, then 
\[
\mathcal{U}\mathfrak{h}_N\cong \left(\mathbb{Q}\langle e_0,e_{\mu_N}\rangle,\circ \right),
\]
where $\circ$ denotes the new product on $\mathbb{Q}\langle e_0,e_{\mu_N}\rangle$ which is transformed from the natural concatenation product on $\mathcal{U}\mathfrak{h}_N$. 

By the same reason as Proposition $2.2$ in \cite{depth}, for any $a\in \mathfrak{h}$, any word sequence $w$ in $e_0,e_{\epsilon},\epsilon\in\mu_N$ and any $n\geq 0$, we have 
\[
a\circ \left(e_0^ne_{\epsilon}w  \right)=e_0^n\left[ \left([\epsilon](a) \right)e_{\epsilon}+e_{\epsilon}\left([\epsilon](a)\right)^*\right]w+e_0^ne_{\epsilon}\left( a\circ w\right),\epsilon\in \mu_N,
\]
where
\[
a\circ e_0^n=e_0^n a,\epsilon\in \mu_N,
\]
\[
\left(u_1 u_2\cdots u_n \right)^*=(-1)^nu_n\cdots u_2 u_1, u_i\in \{e_0,e_{\epsilon};\epsilon\in\mu_N \},
\]
\[
[\epsilon]\left(e_0^{n_1}e_{\epsilon_1}e_0^{n_2}e_{\epsilon_2}\cdots e_0^{n_r}e_{\epsilon_r}e_0^{n_{r+1}} \right)=e_0^{n_1}e_{\epsilon\epsilon_1}e_0^{n_2}e_{\epsilon\epsilon_2}\cdots e_0^{n_r}e_{\epsilon\epsilon_r}e_0^{n_{r+1}},\epsilon,\epsilon_i\in \mu_N.
\]

From the correspondence between unipotent algebraic group and nilpotent Lie algebra (for example, see Section $3$ in \cite{li}), we know that for $a\in \mathfrak{h}_N$, the natural action of $a$ on $\mathcal{O}({}_0\Pi_1)$:
\[
\mathcal{O}({}_0\Pi_1)=\mathbb{Q}\langle e^0,e^{\mu_N}\rangle \xrightarrow{a} \mathcal{O}({}_0\Pi_1)=\mathbb{Q}\langle e^0,e^{\mu_N}\rangle,
\]
\[
x\mapsto a(x),
\]
is dual to the following action of $a$ on $\mathcal{U}\mathfrak{h}$:
\[
\mathcal{U}\mathfrak{h}_N=\mathbb{Q}\langle e_0,e_{\mu_N} \rangle\xrightarrow{a}\mathcal{U}\mathfrak{h}_N=\mathbb{Q}\langle e_0,e_{\mu_N} \rangle,
\]
\[
y\mapsto a\circ y.
\]

By the definition of $\mathcal{H}$ and $\partial_{2n+1}$, we have the following commutative diagram
 \[
 \xymatrix{
   gr_r^{\mathcal{D}}\mathbb{Q}\langle e^0,e^{\mu_N}\rangle \ar@{->>}[d] \ar[r]^{\overline{\partial_{n}}} & gr_{r-1}^{\mathcal{D}}\mathbb{Q}\langle e^0,e^{\mu_N}\rangle \ar@{->>}[d] \\
 gr_r^{\mathcal{D}}\mathcal{H} \ar[r]^{\partial_{n}} & gr_{r-1}^{\mathcal{D}}\mathcal{H} , }
 \]
 where $\overline{\partial_{n}}$ is  the depth-graded version of the action of $i(\sigma_{n})$ on $\mathbb{Q}\langle e^0,e^{\mu_N}\rangle$.

Let $\delta\binom{x}{y}$ be the function of $x,y\in \mathbb{C}$ which satisfies 
\[
\delta\binom{x}{y}=
\begin{cases}
1, &x=y;\\
0, &x\neq y.
\end{cases}
\]
Denote by $\mathfrak{g}_N^{ab}=\mathfrak{g}_N/[\mathfrak{g}_N,\mathfrak{g}_N]$ and $\left(\mathfrak{g}_N^{ab} \right)^{{\vee}}$ be its compact dual (the set of linear functions which are non-zero only for a finite dimension subspace). For $N=2$, let $$f_{2n+1},n\geq 0,\in \left(\mathfrak{g}_N^{ab} \right)^{{\vee}}$$ be the dual basis of the images of $\sigma_{2n+1},n\geq 0$ in $\mathfrak{g}_N^{ab}$. For $N=3,4$, let $$f_{n},n\geq 1,\in \left(\mathfrak{g}_N^{ab} \right)^{{\vee}}$$ be the dual basis of the images of $\sigma_{n},n\geq 1$ in $\mathfrak{g}_N^{ab}$.

For $N=2$, there is a well-defined map
\[
\partial:gr_r^{\mathcal{D}}\mathcal{H}\rightarrow  \left(\mathfrak{g}_2^{ab} \right)^{{\vee}}\otimes gr_{r-1}^{\mathcal{D}}\mathcal{H},\partial=\sum_{n\geq 0}f_{2n+1}\otimes \partial_{2n+1}.
\]
For $N=3,4$, there is a well-defined map
\[
\partial:gr_r^{\mathcal{D}}\mathcal{H}\rightarrow  \left(\mathfrak{g}_N^{ab} \right)^{{\vee}}\otimes gr_{r-1}^{\mathcal{D}}\mathcal{H},\partial=\sum_{n\geq 1}f_{n}\otimes \partial_{n}.
\]

Now we have 
\begin{prop}\label{inj}
For $r\geq 2$, the map $\partial$ is injective for $N=2,3,4$.
\end{prop}
 \noindent{\bf Proof}:
By exactly the same method in Section $2.3$,\cite{brown}, it follows that 
\[
\mathcal{H}\cong \mathcal{O}\left(\mathcal{U}^{\mathcal{MT}_N}\right)[t]
\]
as a $\mathfrak{g}_N$-module, where $t$ is a weight $\begin{cases} 2, &N=2\\1, &N=3,4\end{cases}$, depth $1$ element with trivial action of $\mathfrak{g}_N$.  Furthermore, $t^n,n\geq 1$ are all depth $1$ elements.

As a result,
\[
gr_r^{\mathcal{D}}\mathcal{H}\cong gr_r^{\mathcal{D}}\mathcal{O}\left(\mathcal{U}^{\mathcal{MT}_N}\right)\oplus \bigoplus_{n\geq 1} gr_{r-1}^{\mathcal{D}}\mathcal{O}\left(\mathcal{U}^{\mathcal{MT}_N}\right)t^n.
\]
Be aware that $gr_r^{\mathcal{D}}\mathcal{O}(\mathcal{U}^{\mathcal{MT}_N})$ is dual to $gr_{\mathcal{D}}^r\mathcal{U}\mathfrak{g}_N$ and the decreasing depth filtration on $\mathcal{U}\mathfrak{g}_N$ is induced by the depth filtration on $\mathfrak{g}_N$.

Thus it suffices to prove that $\partial|_{gr_r^{\mathcal{D}}\mathcal{O}\left(\mathcal{U}^{\mathcal{MT}_N}\right)}$ is injective. Since the depth-graded motivic Lie algebra $\mathfrak{dg}$ is a free Lie algebra with generators which are all in the depth one parts \cite{deligne}.
By the correspondence between nilpotent Lie algebra and unipotent algebraic group, $\partial|_{gr_r^{\mathcal{D}}\mathcal{O}\left(\mathcal{U}^{\mathcal{MT}_N}\right)}$ is injective.  
$\hfill\Box$\\

\section{Main results} \label{mainsection}

In this section, firstly we give the set of elements which generate all the cyclotomic multiple zeta values for $\mu_N$ in cases $N=2,3,4$. These elements are $\mathbb{Q}$-linear combinations of unit cyclotomic multiple zeta values for $\mu_N$. Secondly we give the structure of depth-graded unit cyclotomic multiple zeta values. 

\begin{lem}\label{home}
For $N=3,4$, $\pm1\neq \varepsilon\in \mu_N$, the map \[
\gamma: [0,1]\rightarrow \mathbb{C};\;\;\;t\mapsto s=\frac{1-t}{1-\varepsilon t}\]
satisfies:\\
$(i)$ $\gamma(t)\subseteq \mathbb{C}-\{0,\mu_N\},\; \forall\; t\in (0,1)$;\\
$(ii)$ The map $\gamma$ is (fixed base point) homotopically equivalent to the map
\[
\beta:[0,1]\rightarrow\mathbb{C};\;\;\; t\mapsto 1-t
\]
over $\mathbb{C}-\{0,\mu_N\}$.
\end{lem}
\noindent{\bf Proof}: 
Denote by $\varepsilon=cos \, \theta +i\, sin\, \theta $, $\theta \in (0, 2\pi)$. For $t\in (0,1)$, one has 
\[
\begin{split}
&\;\;\;\;\lvert 1-\varepsilon t\rvert=\sqrt{(1-t \,cos\theta)^2+t^2 sin^2\theta}=\sqrt{1-2t\, cos\theta+t^2 }\\
&>\sqrt{1-2t+t^2}=1-t.         \\
\end{split}
\]
Here the above inequality follows from the fact $$cos\,\theta <1, \theta\in (0,2\pi).$$
As a result, $0<\big{\lvert }\frac{1-t}{1-\varepsilon t}\big{\rvert}<1, \forall\, t\in (0,1)$. Thus $(i)$ is proved.

It is easy to check that $u(1-t)+(1-u) \frac{1-t}{1-\varepsilon t}\neq 0, \forall\, t\in (0,1), u\in [0,1]$. So we have
\[
0<\Big{\lvert} u(1-t)+(1-u) \frac{1-t}{1-\varepsilon t}\Big{\rvert}<1, \;\forall\, t\in (0,1), u\in [0,1].
\]
Thus one can define the following map
\[
F: [0,1]\times [0,1]\rightarrow \mathbb{C},
\]
\[
(t,u)\mapsto u(1-t)+(1-u) \frac{1-t}{1-\varepsilon t}.
\]
The map $F$ satisfies that 
$F(0,u)=1,\; F(1,u)=0$ and 
\[
F\big{|}_{u=0}=\gamma, \;F\big{|}_{u=1}=\beta.\]
The statement $(ii)$ is proved.       $\hfill\Box$\\
  
Now we are ready to prove our main results:
\begin{Thm}\label{ubas}
$(i)$ For $N=2$, the $\mathbb{Q}$-linear space of cyclotomic multiple zeta values for $\mu_2$ is generated by the following integrals 
\[
\mathop{\int}_{0<s_1<\cdots<s_r<1}\omega_{i_1}(s_1)\cdots \omega_{i_r}(s_r)\cdot \left(  \mathop{\int}_{0<v_1<v_2<1}   \omega_1(v_1)\omega_2(v_2)\right)^{2m},\]
\[ i_1, \cdots, i_{r-1}\in\{1,2\}, i_r=2,  r\geq 1, m\geq 0.
\]
Here 
\[
\omega_1(s)=\frac{ds}{1+s}+\frac{ds}{1-s},\;\; \omega_2(s)=\frac{ds}{1+s}.
\]\\
$(ii)$ For $N=3, 4$, denote by $\varepsilon $ a primitive root of $N^{th}$-unity.  Then the $\mathbb{Q}$-linear space of cyclotomic multiple zeta values  is generated by the following integrals
\[
\mathop{\int}_{0<s_1<\cdots<s_r<1}\omega_{i_1}(s_1)\cdots \omega_{i_r}(s_r)\cdot \left(  \mathop{\int}_{0}^1   \omega_3(v)\right)^{m},\]
\[
 i_1, \cdots, i_{r-1}\in\{1,2\}, i_r=2,  r\geq 1, m\geq 0.
\]
Here 
\[
\omega_1(s)=\frac{ds}{\varepsilon^{-1}-s}-\frac{ds}{1-s},\;\omega_2(s)=\frac{ds}{\varepsilon^{-1}-s},\; \omega_3(s)=\frac{ds}{\varepsilon^{-1}-s}-\frac{ds}{\varepsilon-s}.\]
\end{Thm}
 \noindent{\bf Proof}: 
 By the Corollary $1.1$ in \cite{gla}, it follows that:\\
 $(A)$ For $N=2$, the set of following elements 
 \[
 \zeta\binom{x_1,\cdots,x_{p-1},x_p}{1\;,\cdots,\,1\;\,,\;-1}\left[\zeta\binom{2}{-1}\right]^m,\; x_i\geq 1\;\mathrm{odd}, m\geq 1.
 \]
 generates all the cyclotomic multiple zeta values for $\mu_2$.\\
 $(B)$ For $N=3,4$, the set of following elements
 \[
 \zeta\binom{x_1,\cdots,x_{p-1},x_p}{1\;,\cdots,\,\,1\,\;\,,\,\;\varepsilon}\left[\zeta\binom{1}{\varepsilon}-\zeta\binom{1}{\varepsilon^{-1}}\right]^m,\; x_i\geq 1, m\geq 1.
 \]
 generates all the cyclotomic multiple zeta values for $\mu_2$.
 
 $(i)$ From the iterated integral representations of cyclotomic multiple zeta values, one has
 \[
  \zeta\binom{x_1,\cdots,x_{p-1},x_p}{\,1\;,\cdots,\,1\;\,,\;\;-1}=\mathop{\int}_{0<t_1<\cdots<t_r<1}\delta_{1}(t_1)\cdots \delta_{r}(t_r). \]
  Here $r=x_1+\cdots +x_p$ and 
  \[
  \delta_{j}(t)=\begin{cases}
   \frac{dt}{-1-t},&j\in \{1, 1+x_1,\cdots, 1+x_1+\cdots+x_{p-1}\}\\
   \frac{dt}{t},&j\notin \{1, 1+x_1,\cdots, 1+x_1+\cdots+x_{p-1}\}.
      \end{cases}
  \]
  By changing of variables 
  \[
  t_1=\frac{1-s_1}{1+s_1}, \cdots, t_r=\frac{1-s_r}{1+s_r},
  \]
  it follows that
   \[
  \delta_{j}(t)=\begin{cases}
   \frac{dt}{-1-t}=\frac{ds}{1+s} ,&j\in \{1, 1+x_1,\cdots, 1+x_1+\cdots+x_{p-1}\}\\
   \frac{dt}{t}=-\frac{ds}{1-s}-\frac{ds}{1+s},&j\notin \{1, 1+x_1,\cdots, 1+x_1+\cdots+x_{p-1}\}.
      \end{cases}
  \]
  As a result, 
  \[
   \zeta\binom{x_1,\cdots,x_{p-1},x_p}{1\;,\cdots,\,1\;\,,\;-1}=\mathop{\int}_{1>s_1>\cdots>s_r>0}u_{1}(s_1)\cdots u_{r}(s_r).
     \]
    Here $r=x_1+\cdots +x_p$ and 
  \[
  u_{j}(s)=\begin{cases}
   \frac{ds}{1+s},&j\in \{1, 1+x_1,\cdots, 1+x_1+\cdots+x_{p-1}\}\\
  - \frac{ds}{1-s}-\frac{ds}{1+s},&j\notin \{1, 1+x_1,\cdots, 1+x_1+\cdots+x_{p-1}\}.
      \end{cases}
  \]
 Thus the statement $(i)$ is proved.
 
 $(ii)$ Similarly, for $N=3,4$, one has
  \[
  \zeta\binom{x_1,\cdots,x_{p-1},x_p}{1\;,\cdots,\,1\;\,,\;\;\,\varepsilon}=\mathop{\int}_{0<t_1<\cdots<t_r<1}\lambda_{1}(t_1)\cdots \lambda_{r}(t_r). \]
  Here $r=x_1+\cdots +x_p$ and 
  \[
  \lambda_{j}(t)=\begin{cases}
   \frac{dt}{\varepsilon^{-1}-t},&j\in \{1, 1+x_1,\cdots, 1+x_1+\cdots+x_{p-1}\}\\
   \frac{dt}{t},&j\notin \{1, 1+x_1,\cdots, 1+x_1+\cdots+x_{p-1}\}.
      \end{cases}
  \]
  By changing of variables 
  \[
  t_1=\frac{1-s_1}{1-\varepsilon s_1}, \cdots, t_r=\frac{1-s_r}{1-\varepsilon s_r},
  \]
   one can check that 
   \[
   s_1=\frac{1-t_1}{1-\varepsilon t_1}, \cdots, s_r=\frac{1-t_r}{1-\varepsilon t_r}.   \]
   It is easy to check that 
   \[
   \frac{dt}{\varepsilon^{-1}-t}=-\frac{ds}{\varepsilon^{-1}-s},\;\frac{dt}{t}=\frac{ds}{\varepsilon^{-1}-s}-\frac{ds}{1-s}.
   \]
   So we have 
   \[
  \zeta\binom{x_1,\cdots,x_{p-1},x_p}{1\;,\cdots,\,1\;\,,\;\;\,\varepsilon}=\mathop{\int\cdots \int}_{\Delta_r}v_{1}(s_1)\cdots v_{r}(s_r). \]
  Here $r=x_1+\cdots +x_p,$
  \[
  \begin{split}
  &\;\;\;\;\Delta_r\\
  &=\Bigg{\{}\left(\frac{1-t_1}{1-\varepsilon t_1},\cdots, \frac{1-t_r}{1-\varepsilon t_r} \right)\in \mathbb{C}^r\, \Big{|}\,0<t_1<\cdots<t_r<1     \Bigg{\}}        \\
  &=  \Bigg{\{}\left(s_1,\cdots,s_r \right)\in \mathbb{C}^r \,\Big{|}\,\frac{1-s_i}{1-\varepsilon s_i}\in \mathbb{R},\forall\, i,\, 0< \frac{1-s_1}{1-\varepsilon s_1}<\cdots< \frac{1-s_r}{1-\varepsilon s_r} <1  \Bigg{\}}             \\
  \end{split}
  \]
   and 
  \[
  v_{j}(s)=\begin{cases}
   -\frac{ds}{\varepsilon^{-1}-s},&j\in \{1, 1+x_1,\cdots, 1+x_1+\cdots+x_{p-1}\},\\
       \frac{ds}{\varepsilon^{-1}-s}-\frac{ds}{1-s}   ,&j\notin \{1, 1+x_1,\cdots, 1+x_1+\cdots+x_{p-1}\}.
      \end{cases}
  \]
  By Lemma \ref{home}, $\Delta_r$ is homotopically equivalent to 
  \[
  \begin{split}
  &\;\;\;\; \Gamma_r\\
  &=\Big{\{}\left(1-t_1,\cdots, 1-t_r \right)\in \big{(}\mathbb{C}-\{0,\mu_N\}\big{)}^r\, \Big{|}\,0<t_1<\cdots<t_r<1     \Big{\}}        \\
  &=  \Big{\{}\left(s_1,\cdots,s_r \right)\in \big{(}\mathbb{C}-\{0,\mu_N\}\big{)}^r  \,\Big{|}\, 1>s_1>\cdots>s_r>0 \Big{\}}             \\
  \end{split}
  \]
  over $\big{(}\mathbb{C}-\{0,\mu_N\}\big{)}^r$. Since the differential $v_{1}(s_1)\cdots v_{r}(s_r)$ is holomorphic over $\big{(}\mathbb{C}-\{0,\mu_N\}\big{)}^r$, we have
  \[
    \zeta\binom{x_1,\cdots,x_{p-1},x_p}{1\;,\cdots,\,1\;\,,\;\;\,\varepsilon}=\mathop{\int}_{1>s_1>\cdots>s_r>0}v_{1}(s_1)\cdots v_{r}(s_r).   \]   
Thus the statement $(i)$ is proved.   $\hfill\Box$\\
 
 \begin{rem}
 In the above proof, the identity
 \[
 \mathop{\int\cdots \int}_{\Delta_r}v_{1}(s_1)\cdots v_{r}(s_r)-\mathop{\int\cdots \int}_{\Gamma_r} v_{1}(s_1)\cdots v_{r}(s_r)  =0\]
 can be viewed as a special case of the following generalized Stokes' formula on manifold:
 \[
\int_{\partial M}\omega= \int_Md\omega.
 \]
 Here 
 \[
 \begin{split}
  &M=\Big{\{}(s_1,\cdots, s_r,u)\in \mathbb{C}^{r}\times (0,1)\,\Big{|}\,\substack{ s_i=\frac{1-t_i}{1-\varepsilon t_i}u+(1-t_i)(1-u), \forall \,i\\0<t_1<\cdots<t_r<1,\; u\in (0,1)} \Big{\}},\\
 \end{split}
 \]
 \[
 \omega=v_1(s_1)\cdots v_r(s_r).
 \]
 \end{rem}
 \begin{rem}\label{tone}
 From the proof of Theorem \ref{ubas}, one has 
 \[
 \zeta\binom{2}{-1}=-\zeta\binom{1\,,\;1}{1,-1}+\zeta\binom{1\;,\;1}{-1,-1}.
 \]
 \end{rem}
 
 By the shuffle product of cyclotomic multiple zeta values, Theorem \ref{uzhao} follows immediately from Theorem \ref{ubas}.
 
\begin{Thm}\label{mmm}
$(i)$ For $N=2,r\geq 1$, $\mathrm{dim}_{\mathbb{Q}}\,gr_r^{\mathcal{D}}\mathcal{H}_r^{(1)}=1$ and $gr_r^{\mathcal{D}}\mathcal{H}^{(1)}_r$ is generated by 
\[
\zeta^{\mathfrak{m}}\dbinom{{1,\cdots,1,\;\;1}}{\underbrace{1,\cdots,1,-1}_r}
\]
as a $\mathbb{Q}$-linear space;\\
$(ii)$ For $N=3,4,r\geq 1$, $\mathrm{dim}_{\mathbb{Q}}\,gr_r^{\mathcal{D}}\mathcal{H}_r^{(1)}=2$ and $gr_r^{\mathcal{D}}\mathcal{H}^{(1)}_r$ is generated by \[
\zeta^{\mathfrak{m}}\dbinom{{1,\cdots,1,1}}{\underbrace{1,\cdots,1,\epsilon}_r},\zeta^{\mathfrak{m}}\dbinom{{1,\cdots,1,\;1\;\;}}{\underbrace{1,\cdots,1,\epsilon^{-1}}_r}
\]
 as a $\mathbb{Q}$-linear space.
\end{Thm}
 \noindent{\bf Proof}: For $r=1$, it is clear that (i) and (ii) are true by definition.
 Since the map $\partial$ is injective, from Proposition \ref{inj} and Lemma \ref{below} below, it follows that $\partial_1$ is injective for $\mu_2,\mu_3$ and $\mu_4$. Thus we have
 \[
 \mathrm{dim}_{\mathbb{Q}}gr_r^{\mathcal{D}}\mathcal{H}_r^{(1)}=\mathrm{dim}_{\mathbb{Q}}\underbrace{\partial_1\circ \cdots\circ\partial_1}_{r-1}\left( gr_r^{\mathcal{D}}\mathcal{H}^{(1)}\right).
 \]
 From the explicit formulas of $\partial_1$ in Lemma \ref{below}, we have for $N=2$,
 \[
 \underbrace{ \partial_1\circ\cdots\circ\partial_1}_{r-1}\left( \zeta^{\mathfrak{m}}\dbinom{{1,\cdots,1,\,\;1\;}}{\underbrace{1,\cdots,1,-1}_r} \right)=\zeta^{\mathfrak{m}}\binom{1}{-1}
  \]
  and for $N=3$ and $\epsilon $ an primitive of N-th unity  \[
   \underbrace{ \partial_1\circ\cdots\circ\partial_1}_{r-1}\left(\zeta^{\mathfrak{m}}\dbinom{{1,\cdots,1,\;1\;\;}}{\underbrace{1,\cdots,1,\epsilon^{\pm1}}_r}
\right)=\zeta^{\mathfrak{m}}\binom{1}{\epsilon^{\pm1}}.
     \]
 Thus for $N=2,3$ the theorem is proved.  
 For $N=4$ and $\epsilon $ an primitive of N-th unity, 
  \[
   \underbrace{ \partial_1\circ\cdots\circ\partial_1}_{r-1}\left(\zeta^{\mathfrak{m}}\dbinom{{1,\cdots,1,\;1\;\;}}{\underbrace{1,\cdots,1,\epsilon^{\pm1}}_r}
\right)=\zeta^{\mathfrak{m}}\binom{1}{\epsilon^{\pm1}}.
     \]
 and 
 \[
 \underbrace{ \partial_1\circ\cdots\circ\partial_1}_{r-1}\left( \zeta^{\mathfrak{m}}\dbinom{{1,\cdots,1,\;\;1}}{\underbrace{1,\cdots,1,-1}_r} \right)=2^{r-1}\zeta^{\mathfrak{m}}\binom{1}{-1}
  \]
 For $N=4$, the theorem follows from 
 \[
 \zeta^{\mathfrak{m}}\binom{1}{\epsilon}+\zeta^{\mathfrak{m}}\binom{1}{\epsilon^{-1}}=\zeta^{\mathfrak{m}}\binom{1}{-1}.
 \]  $\hfill\Box$\\
 
 \begin{rem}
 For $N=4$, $\mu_4=\{\pm1, \pm i\}$, one has 
 \[
 \zeta\binom{1}{i}+\zeta\binom{1}{-i}=\int^1_0\left( \frac{1}{-i-t}+\frac{1}{i-t}  \right)dt=\int^1_0-\frac{2t}{1+t^2}dt=\zeta\binom{1}{-1}.
 \]
 \end{rem}

From Theorem \ref{mmm}, by the period map $per:\mathcal{H}\rightarrow \mathbb{C}$ we get
 Theorem \ref{gene} immediately.

\begin{lem}\label{below}
(i) For $N=2,n\geq 1$, $\partial_{2n+1}\left(gr_1^{\mathcal{D}}\mathcal{H}^{(1)} \right)=0$. For $e^{i_1},\cdots, e^{i_s}\in \{\pm1\}$, we have
\[
\begin{split}
&\;\;\;\;\overline{\partial_1}\left(e^{i_1}e^{i_2}\cdots e^{i_s} \right)\\
&=\delta\binom{i_1}{-i_2}e^{i_2}\cdots e^{i_s} +\cdots +\delta\binom{i_{s-1}}{-i_s}e^{i_1}\cdots e^{i_{s-2}}e^{i_s}+\delta\binom{i_s}{-1}e^{i_1}\cdots e^{i_{s-1}}   \\
&\;\;\;\;-\delta\binom{i_1}{-i_2}e^{i_1}e^{i_3}\cdots e^{i_s}-\cdots-\delta\binom{i_{s-1}}{-i_s}e^{i_1}\cdots e^{i_{s-1}}  \\
\end{split}
\]
(ii) For $N=3,n\geq 2$, $\partial_{n}\left(gr_1^{\mathcal{D}}\mathcal{H}^{(1)}\right)=0$.  For $e^{i_1},\cdots, e^{i_s}\in \mu_3$, we have
\[
\begin{split}
&\;\;\;\;\overline{\partial_1}\left(e^{i_1}e^{i_2}\cdots e^{i_s}\right)\\
&=\left[ \delta\binom{i_1}{i_2\epsilon}+\delta\binom{i_1}{i_2\epsilon^{-1}}\right]e^{i_2}\cdots e^{i_s}+\cdots+\left[\delta\binom{i_{s-1}}{i_s\epsilon}+\delta\binom{i_{s-1}}{i_s\epsilon^{-1}} \right]e^{i_1}\cdots e^{i_{s-2}}e^{i_s}    \\
&+\left[\delta\binom{i_s}{\epsilon}+\delta\binom{i_s}{\epsilon^{-1}} \right]e^{i_1}\cdots e^{i_{s-1}}-\left[\delta\binom{i_1}{i_2\epsilon}+\delta\binom{i_1}{i_2\epsilon^{-1}}\right]e^{i_1}e^{i_3}\cdots e^{i_s}\\
&\cdots- \left[\delta\binom{i_{s-1}}{i_s\epsilon}+\delta\binom{i_{s-1}}{i_s\epsilon^{-1}} \right]e^{i_1}\cdots e^{i_{s-2}} e^{i_{s-1}} .  \\
\end{split}
\]
(iii) For $N=4, n\geq 2$, $\partial_{n}\left(gr_1^{\mathcal{D}}\mathcal{H}^{(1)}\right)=0$. For $e^{i_1},\cdots, e^{i_s}\in \mu_4$, we have
\[\tiny
\begin{split}
&\;\;\;\;\overline{\partial_1}\left(e^{i_1}e^{i_2}\cdots e^{i_s}\right)\\
&=\left[ \delta\binom{i_1}{i_2\epsilon}+2\delta\binom{i_1}{-i_2}+\delta\binom{i_1}{i_2\epsilon^{-1}}\right]e^{i_2}\cdots e^{i_s}+\cdots+\left[\delta\binom{i_{s-1}}{i_s\epsilon}+2\delta\binom{i_{s-1}}{-i_s}+\delta\binom{i_{s-1}}{i_s\epsilon^{-1}} \right]e^{i_1}\cdots e^{i_{s-2}}e^{i_s}    \\
&+\left[\delta\binom{i_s}{\epsilon}+2\delta\binom{i_s}{-1}+\delta\binom{i_s}{\epsilon^{-1}} \right]e^{i_1}\cdots e^{i_{s-1}}-\left[\delta\binom{i_1}{i_2\epsilon}+2\delta\binom{i_1}{-i_2}+\delta\binom{i_1}{i_2\epsilon^{-1}} \right]e^{i_1}e^{i_3}\cdots e^{i_s}\\
&-\cdots- \left[\delta\binom{i_{s-1}}{i_s\epsilon}+2\delta\binom{i_{s-1}}{-i_s}+\delta\binom{i_{s-1}}{i_s\epsilon^{-1}} \right]e^{i_1}\cdots e^{i_{s-2}}e^{i_{s-1}}.  \\
\end{split}
\]
\end{lem}
 \noindent{\bf Proof}:
 (i) From the commutative diagram in Section \ref{mga}, to prove that $$\partial_{2n+1}\left(gr_1^{\mathcal{D}}\mathcal{H}^{(1)} \right)=0, \forall\, n\geq 1$$
 it suffices to prove that 
\[
\overline{\partial_{2n+1}}\left(\mathbb{Q}\langle e^{\mu_N}\rangle\right)=0.
\]
Here $\mathbb{Q}\langle e^{\mu_n}\rangle$ is the sub-algebra of $\mathbb{Q}\langle e^0,e^{\mu_N}\rangle$ generated by $e^{\epsilon_1}e^{\epsilon_2}\cdots e^{\epsilon_r},\epsilon_i\in \mu_N,r\geq 1$.
By considering the action of $\overline{\sigma_{2n+1}}$ on $\mathcal{U}\mathfrak{h}=\mathbb{Q}\langle e_0,e_{\mu_N}\rangle$, from Section \ref{mga}, it is enough to show that the terms
\[
e_{\xi_1}e_{\xi_2}\cdots e_{\xi_{r+1}}, \xi_1,\cdots,\xi_{r+1}\in \mu_N
\]
have trivial coefficients in 
\[
\overline{\sigma_{2n+1}}\circ e_{\epsilon_1}e_{\epsilon_2}\cdots e_{\epsilon_r}, \forall\, \epsilon_1,\cdots,\epsilon_r\in \mu_N
\]
for all $r\geq 0$. This follows from the definition of $\circ$ and $\overline{\sigma_{2n+1}}$. While the formula for $\overline{\partial_1}$ follows from that
\[
\begin{split}
&\;\;\;\;e_{-1}\circ\left(e_{i_1}e_{i_2}\cdots e_{i_r}\right)\\
&=\left(e_{-i_1}e_{i_1}-e_{i_1}e_{-i_1} \right)e_{i_2}\cdots e_{i_r}+e_{i_1}\left( e_{-i_2}e_{i_2}-e_{i_2}e_{-i_2}\right)e_{i_3}\cdots e_{i_r}+\cdots \\
&+e_{i_1}\cdots e_{i_{r-1}}\left(e_{-i_r}e_{i_r} -e_{i_r}e_{-i_r}\right)+e_{i_1}\cdots e_{i_r}e_{-1}.\\
\end{split}
\]

Similarly, for $N=3,4$, $\eta\in \mu_N$ and $\eta\neq 1$, the statements $(ii)$ and $(iii)$ follow from 
\[
\begin{split}
&\;\;\;\;e_{\eta}\circ\left(e_{i_1}e_{i_2}\cdots e_{i_r}\right)\\
&=\left(e_{\eta i_1}e_{i_1}-e_{i_1}e_{\eta i_1} \right)e_{i_2}\cdots e_{i_r}+e_{i_1}\left( e_{\eta i_2}e_{i_2}-e_{i_2}e_{\eta i_2}\right)e_{i_3}\cdots e_{i_r}+\cdots \\
&+e_{i_1}\cdots e_{i_{r-1}}\left(e_{\eta i_r}e_{i_r} -e_{i_r}e_{\eta i_r}\right)+e_{i_1}\cdots e_{i_r}e_{\eta }\\ 
\end{split}\]
 and the the expression of $\overline{\partial_{2n+1}}$ in cases $N=3,4$.
 $\hfill\Box$\\
 
 \begin{rem}
 For general $N$, is the $\mathbb{Q}$-algebra of $\bf unit$ cyclotomic multiple zeta values for $\mu_N$  still equal to the $\mathbb{Q}$-algebra of  cyclotomic multiple zeta values for $\mu_N$?
  \end{rem}
  \begin{rem}
  For any $N$, all the weight one unit cyclotomic multiple zeta values are known to be transcendental. Up to now, we nearly know nothing about the  unit cyclotomic multiple zeta values of weight $\geq 2$.
  \end{rem}

\section{Calculations of unit cyclotomic multiple zeta values}
In this section, we calculate the coefficients of the unit cyclotomic multiple zeta values in terms of depth-graded basis. Furthermore, we give some examples about  the explicit relations among (depth-graded) unit cyclotomic multiple zeta values.

Since 
\[
\begin{split}
&\zeta\dbinom{1,\cdots,1,1}{\underbrace{1,\cdots,1,\epsilon}_r}
= dch(\underbrace{e^{\epsilon^{-1}}e^{\epsilon^{-1}}\cdots e^{\epsilon^{-1}}}_r ) 
=  \frac{1}{r!}dch(\underbrace{e^{\epsilon^{-1}}\rotatebox{90}{$\rotatebox{180}{$\exists$}$}\,e^{\epsilon^{-1}}\rotatebox{90}{$\rotatebox{180}{$\exists$}$}\,\cdots \rotatebox{90}{$\rotatebox{180}{$\exists$}$}\,e^{\epsilon^{-1}}}_r ) 
=\frac{1}{r!}\left( \zeta\binom{1}{\epsilon} \right)^r,\\
\end{split}
\]
we have 
\[
\zeta\dbinom{1,\cdots,1,1}{\underbrace{1,\cdots,1,\epsilon}_r}=\frac{(-1)^r}{r!}\left[\mathrm{log}\;(1-\epsilon) \right]^r.
\]
From Lemma \ref{below}, for any 
\[
\zeta\binom{1,\;\cdots,1}{\epsilon_1,\cdots,\epsilon_r}
\]
one can use the formulas for $\overline{\partial_1}$ inductively to calculate the numbers $$a_{\epsilon_1,\cdots,\epsilon_r},b_{\epsilon_1,\cdots,\epsilon_r}, c_{\epsilon_1,\cdots,\epsilon_r}$$ in the introduction. The key idea is to calculate the map $\underbrace{\overline{\partial_1}\circ\cdots \circ \overline{\partial_1}}_{r-1}$ explicitly for any $r\geq 2$.

\begin{Def}
For $N=2,3,4$, define the $\mathbb{Q}$-linear maps $$\mathcal{D}_1,\mathcal{D}_2: \mathbb{Q}\langle e^{\mu_N}\rangle \rightarrow  \mathbb{Q}\langle e^{\mu_N}\rangle$$ as 
\[
\mathcal{D}_1(1)=0,\;\mathcal{D}_1(e^{i_1})=\delta\binom{i_1}{1},
\]
\[
\begin{split}
&\;\;\;\;\mathcal{D}_1(e^{i_1}e^{i_2}\cdots e^{i_r})\\
&=\delta\binom{i_1}{i_2}e^{i_2}e^{i_3}\cdots e^{i_r}+\cdots+\delta\binom{i_{r-1}}{i_r} e^{i_1}\cdots e^{i_{r-2}}e^{i_r}+\delta\binom{i_r}{1} e^{i_1}\cdots e^{i_{r-2}}e^{i_{r-1}}\\
&\;\;\;-\delta\binom{i_1}{i_2}e^{i_1}e^{i_3}\cdots e^{i_r}-\cdots-\delta\binom{i_{r-1}}{i_r} e^{i_1}\cdots e^{i_{r-2}}e^{i_{r-1}},\end{split}
\]
\[
\mathcal{D}_2(1)=0,\;\mathcal{D}_2(e^{i_1})=1,\;\mathcal{D}_2(e^{i_1}e^{i_2}\cdots e^{i_r})=e^{i_2}\cdots e^{i_r}.
\]
For $N=4$, define the $\mathbb{Q}$-linear map $\mathcal{P}: \mathbb{Q}\langle e^{\mu_N}\rangle \rightarrow  \mathbb{Q}\langle e^{\mu_N}\rangle$ as 
\[
\mathcal{P}(1)=0,\;\mathcal{P}(e^{i_1})=\delta\binom{i_1}{-1},
\]
\[
\begin{split}
&\;\;\;\;\mathcal{P}(e^{i_1}e^{i_2}\cdots e^{i_r})\\
&=\delta\binom{i_1}{-i_2}e^{i_2}e^{i_3}\cdots e^{i_r}+\cdots+\delta\binom{i_{r-1}}{-i_r} e^{i_1}\cdots e^{i_{r-2}}e^{i_r}+\delta\binom{i_r}{-1} e^{i_1}\cdots e^{i_{r-2}}e^{i_{r-1}}\\
&\;\;\;-\delta\binom{i_1}{-i_2}e^{i_1}e^{i_3}\cdots e^{i_r}-\cdots-\delta\binom{i_{r-1}}{-i_r} e^{i_1}\cdots e^{i_{r-2}}e^{i_{r-1}}.\end{split}
\]
\end{Def}

\begin{prop}\label{pd} For convenience, we write the restriction map $$\overline{\partial_1}\Big{|}_{ \mathbb{Q}\langle e^{\mu_N}\rangle}: \mathbb{Q}\langle e^{\mu_N}\rangle\rightarrow \mathbb{Q}\langle e^{\mu_N}\rangle$$
as $\overline{\partial_1}$. \\
$(i)$ For $r\geq 2$, $\mathcal{D}_1(e^{i_1}e^{i_2}\cdots e^{i_r})=\delta\binom{i_r}{1} e^{i_1}\cdots e^{i_{r-1}}$;\\
$(ii)$  For $N=2,3,4$, $\mathcal{D}_1\mathcal{D}_2= \mathcal{D}_2\mathcal{D}_1$;         \\
$(iii)$  For $N=4$,   $\mathcal{D}_1\mathcal{P}=\mathcal{P}\mathcal{D}_1$ and $\mathcal{D}_2\mathcal{P}=\mathcal{P}\mathcal{D}_2$;  \\
$(iv)$ For $N=2,3$, $\overline{\partial_1}+\mathcal{D}_1=\mathcal{D}_2$;\\
$(v)$  For $N=4$,   $\overline{\partial_1}+\mathcal{D}_1=\mathcal{D}_2+\mathcal{P}$;\\
$(vi)$ For $N=4$, $r\geq 1$, $\eta, i_1,\cdots, i_r\in \mu_4$, $\eta/i_1\in \{\epsilon^{\pm 1}\}$,
\[
\mathcal{P}(e^{\eta}e^{i_1}\cdots e^{i_r})=e^{\eta}\mathcal{P}(e^{i_1}\cdots e^{i_r}).\]       
\end{prop}
 \noindent{\bf Proof}:
 $(i)$ The statement follows immediately from the following simple observation
 \[
 \delta\binom{i_1}{i_2}e^{i_2}= \delta\binom{i_1}{i_2}e^{i_1},\cdots,  \delta\binom{i_{r-1}}{i_r}e^{i_r}= \delta\binom{i_{r-1}}{i_r}e^{i_{r-1}}. \]
 $(ii)$ By definition, it is easy to check that 
 \[
 \mathcal{D}_1\mathcal{D}_2(1)=\mathcal{D}_2\mathcal{D}_1(1)=0,
  \]
   \[
 \mathcal{D}_1\mathcal{D}_2(e^i)=\mathcal{D}_2\mathcal{D}_1(e^i)=0, i \in \mu_N.
  \]
 By $(i)$, for $r\geq 2$, one has 
 \[
 \mathcal{D}_1\mathcal{D}_2(e^{i_1}e^{i_2}\cdots e^{i_r})=\mathcal{D}_1(e^{i_2}\cdots e^{i_r})=
 \begin{cases} \delta\binom{i_2}{1}, & r=2,\\
 \delta\binom{i_r}{1}e^{i_2}\cdots e^{i_{r-1}}, & r\geq 3.
 \end{cases}
 \]
 \[
 \mathcal{D}_2\mathcal{D}_1(e^{i_1}e^{i_2}\cdots e^{i_r})=\delta\binom{i_r}{1}\mathcal{D}_2(e^{i_1}\cdots e^{i_{r-1}})=
 \begin{cases} \delta\binom{i_2}{1}, & r=2,\\
 \delta\binom{i_r}{1}e^{i_2}\cdots e^{i_{r-1}}, & r\geq 3.
 \end{cases}
 \]
 As a result, $\mathcal{D}_1\mathcal{D}_2= \mathcal{D}_2\mathcal{D}_1$.\\
 $(iii)$ By definition, one can check that 
  \[
\mathcal{D}_1\mathcal{P}(1)=\mathcal{P}\mathcal{D}_1(1)=\mathcal{D}_2\mathcal{P}(1)=\mathcal{P}\mathcal{D}_2(1) =0,
  \]
   \[
\mathcal{D}_1\mathcal{P}(e^\eta)=\mathcal{P}\mathcal{D}_1(e^\eta)= \mathcal{D}_2\mathcal{P}(e^\eta)=\mathcal{P}\mathcal{D}_2(e^\eta) =0, \eta \in \mu_4,
  \]
    \[
  \begin{split}
  & \;\;\;\;  \mathcal{D}_1\mathcal{P}(e^{i_1}e^{i_2})       \\
  & = \delta\binom{i_1}{-i_2} \mathcal{D}_1(e^{i_2}) +\delta\binom{i_2}{-1}\mathcal{D}_1(e^{i_1})-\delta\binom{i_1}{-i_2}\mathcal{D}_1(e^{i_1})   \\
  &= \delta\binom{i_1}{-i_2} \delta\binom{i_2}{1} +\delta\binom{i_2}{-1}\delta\binom{i_1}{1}-\delta\binom{i_1}{-i_2}\delta\binom{i_1}{1}
  \\
  &=\delta\binom{i_1}{-1}\delta\binom{i_2}{1}\\
  &=\delta\binom{i_2}{1}\mathcal{P}(e^{i_1})\\
  &= \mathcal{P}\mathcal{D}_1(e^{i_1}e^{i_2}),  \;\; i_1,i_2\in \mu_4, \end{split}
  \]
  \[
  \begin{split}
  & \;\;\;\;  \mathcal{D}_2\mathcal{P}(e^{i_1}e^{i_2})       \\
  & = \delta\binom{i_1}{-i_2} \mathcal{D}_2(e^{i_2}) +\delta\binom{i_2}{-1}\mathcal{D}_2(e^{i_1})-\delta\binom{i_1}{-i_2}\mathcal{D}_2(e^{i_1})   \\
  &= \delta\binom{i_1}{-i_2}  +\delta\binom{i_2}{-1}-\delta\binom{i_1}{-i_2}
  \\
  &=\mathcal{P}(e^{i_2})\\
  &= \mathcal{P}\mathcal{D}_2(e^{i_1}e^{i_2}),  \;\; i_1,i_2\in \mu_4. \end{split}
  \]

   By $(i)$, for $r\geq 3$, one has 
     \[\small
   \begin{split}
   &\;\;\;\;\mathcal{D}_1\mathcal{P}(e^{i_1}e^{i_2}\cdots e^{i_r})\\
   &=  \delta\binom{i_1}{-i_2}\mathcal{D}_1(e^{i_2}e^{i_3}\cdots e^{i_r})+\cdots+\delta\binom{i_{r-1}}{-i_r} \mathcal{D}_1(e^{i_1}\cdots e^{i_{r-2}}e^{i_r})+\delta\binom{i_r}{-1} \mathcal{D}_1(e^{i_1}\cdots e^{i_{r-2}}e^{i_{r-1}})\\
&\;\;\;-\delta\binom{i_1}{-i_2}\mathcal{D}_1(e^{i_1}e^{i_3}\cdots e^{i_r})-\cdots-\delta\binom{i_{r-1}}{-i_r} \mathcal{D}_1(e^{i_1}\cdots e^{i_{r-2}}e^{i_{r-1}} )  \\
&=  \delta\binom{i_1}{-i_2}\delta\binom{i_r}{1}e^{i_2}e^{i_3}\cdots e^{i_{r-1}}+\cdots+\delta\binom{i_{r-1}}{-i_r}\delta\binom{i_r}{1} e^{i_1}\cdots e^{i_{r-2}}+\delta\binom{i_r}{-1} \delta\binom{i_{r-1}}{1}e^{i_1}\cdots e^{i_{r-2}}\\
&\;\;\;-\delta\binom{i_1}{-i_2}\delta\binom{i_r}{1}e^{i_1}e^{i_3}\cdots e^{i_{r-1}}-\cdots-\delta\binom{i_{r-1}}{-i_r} \delta\binom{i_{r-1}}{1}e^{i_1}\cdots e^{i_{r-2}}  \\
&=\delta\binom{i_r}{1}\left[  \delta\binom{i_1}{-i_2}e^{i_2}e^{i_3}\cdots e^{i_{r-1}} +\cdots +\delta\binom{i_{r-2}}{-i_{r-1}}e^{i_1}\cdots e^{i_{r-3}}e^{i_{r-1}}+\delta\binom{i_{r-1}}{-1} e^{i_1}\cdots e^{i_{r-2}}     \right]\\
&\;\;\;\;- \delta\binom{i_r}{1}\left[\delta\binom{i_1}{-i_2}e^{i_1}e^{i_3}\cdots e^{i_{r-1}} +\cdots +\delta\binom{i_{r-2}}{-i_{r-1}}e^{i_1}\cdots e^{i_{r-3}}e^{i_{r-2}}        \right]\\
&= \delta\binom{i_r}{1}\mathcal{P}(e^{i_1}e^{i_2}\cdots e^{i_{r-1}})\\
&= \mathcal{P}\mathcal{D}_1(e^{i_1}e^{i_2}\cdots e^{i_{r}}),
\end{split}
\]
   \[\small
   \begin{split}
   &\;\;\;\;\mathcal{D}_2\mathcal{P}(e^{i_1}e^{i_2}\cdots e^{i_r})\\
   &=  \delta\binom{i_1}{-i_2}\mathcal{D}_2(e^{i_2}e^{i_3}\cdots e^{i_r})+\cdots+\delta\binom{i_{r-1}}{-i_r} \mathcal{D}_2(e^{i_1}\cdots e^{i_{r-2}}e^{i_r})+\delta\binom{i_r}{-1} \mathcal{D}_2(e^{i_1}\cdots e^{i_{r-2}}e^{i_{r-1}})\\
&\;\;\;-\delta\binom{i_1}{-i_2}\mathcal{D}_2(e^{i_1}e^{i_3}\cdots e^{i_r})-\cdots-\delta\binom{i_{r-1}}{-i_r} \mathcal{D}_2(e^{i_1}\cdots e^{i_{r-2}}e^{i_{r-1}} )  \\
&=  \delta\binom{i_1}{-i_2}e^{i_3}\cdots e^{i_{r}}+\cdots+\delta\binom{i_{r-1}}{-i_r} e^{i_2}\cdots e^{i_{r-2}}e^{i_r}+\delta\binom{i_r}{-1} e^{i_2}\cdots e^{i_{r-2}}e^{i_{r-1}}\\
&\;\;\;-\delta\binom{i_1}{-i_2}e^{i_3}\cdots e^{i_{r}}-\cdots-\delta\binom{i_{r-1}}{-i_r} e^{i_2}\cdots e^{i_{r-2}}e^{i_{r-1}}  \\
&=  \delta\binom{i_2}{-i_3}e^{i_3}e^{i_4}\cdots e^{i_{r}}+\cdots+\delta\binom{i_{r-1}}{-i_r} e^{i_2}\cdots e^{i_{r-2}}e^{i_r}+\delta\binom{i_r}{-1} e^{i_2}\cdots e^{i_{r-2}}e^{i_{r-1}}\\
&\;\;\;-\delta\binom{i_2}{-i_3}e^{i_2}e^{i_4}\cdots e^{i_{r}}-\cdots-\delta\binom{i_{r-1}}{-i_r} e^{i_2}\cdots e^{i_{r-2}}e^{i_{r-1}}  \\
&=\mathcal{P}(e^{i_2}e^{i_3}\cdots e^{i_r})\\
&=\mathcal{P}\mathcal{D}_2(e^{i_1}e^{i_2}\cdots e^{i_r}) .  \end{split}
   \]
   In a word, one has  $\mathcal{D}_1\mathcal{P}=\mathcal{P}\mathcal{D}_1$ and $\mathcal{D}_2\mathcal{P}=\mathcal{P}\mathcal{D}_2$.\\
$(iv)$ For $N=2$, it is clear that 
\[
\delta\binom{\eta_1}{\eta_2}+\delta\binom{\eta_1}{-\eta_2}=1, \eta_1,\eta_2\in \{\pm 1\}.\]
Thus 
\[
(\overline{\partial_1}+\mathcal{D}_1)(1)=0=\mathcal{D}_2(1),
\]
\[
(\overline{\partial_1}+\mathcal{D}_1)(e^{\eta})=\delta\binom{\eta}{1}+\delta\binom{\eta}{-1} =1=\mathcal{D}_2(e^\eta),\;\eta\in \{\pm 1\}.
\]
For $r\geq 2$, 
\[
\begin{split}
&\;\;\;\; (\overline{\partial_1}+\mathcal{D}_1)(e^{i_1}e^{i_2}\cdots e^{i_r})      \\
&=\left[\delta\binom{i_1}{-i_2}+ \delta\binom{i_1}{i_2}\right]e^{i_2}e^{i_3}\cdots e^{i_r}+\cdots+\left[ \delta\binom{i_{r-1}}{-i_r}+ \delta\binom{i_{r-1}}{i_r}\right] e^{i_1}\cdots e^{i_{r-2}}e^{i_r}\\
&\;\;\;\;+\left[\delta\binom{i_r}{-1}+ \delta\binom{i_r}{1}\right] e^{i_1}\cdots e^{i_{r-2}}e^{i_{r-1}}\\
&\;\;\;\;-\left[\delta\binom{i_1}{-i_2}+ \delta\binom{i_1}{i_2}\right] e^{i_1}e^{i_3}\cdots e^{i_r}-\cdots- \left[ \delta\binom{i_{r-1}}{-i_r}+ \delta\binom{i_{r-1}}{i_r}\right] e^{i_1}\cdots e^{i_{r-2}}e^{i_{r-1}}   \\
&=e^{i_2}e^{i_3}\cdots e^{i_r}+e^{i_1}e^{i_3}\cdots e^{i_r}+\cdots +e^{i_1}\cdots e^{i_{r-2}}e^{i_{r-1}}\\
&\;\;\;\;-e^{i_1}e^{i_3}\cdots e^{i_r}-\cdots -e^{i_1}\cdots e^{i_{r-2}}e^{i_{r-1}}  \\
&=e^{i_2}e^{i_3}\cdots e^{i_r}\\
&=\mathcal{D}_2(   e^{i_1}e^{i_2}\cdots e^{i_r}  ).
\end{split}
\]
By the exactly same analysis, the formula $\overline{\partial_1}+\mathcal{D}_1=\mathcal{D}_2$ in case $N=3$ follows from 
\[
\delta\binom{\eta_1}{\eta_2}+\delta\binom{\eta_1}{\epsilon \eta_2}+\delta\binom{\eta_1}{\epsilon^{-1} \eta_2}=1,\; \eta_1,\eta_2\in \mu_3.\]
$(v)$ For $N=4$, one can check that 
\[
( \overline{\partial_1}+\mathcal{D}_1)(1)=0=(\mathcal{D}_2+\mathcal{P})(1),
 \]
 \[
 \begin{split}
&\;\;\;\; ( \overline{\partial_1}+\mathcal{D}_1)(e^{\eta})\\
&=\delta\binom{\eta}{\epsilon}+2\delta\binom{\eta}{-1}+\delta\binom{\eta}{\epsilon^{-1}} +\delta\binom{\eta}{1}\\
&=1+\delta\binom{\eta}{-1}\\
&=(\mathcal{D}_2+\mathcal{P})(e^{\eta}).\\
\end{split} \]
 For $r\geq 2$,
  \[\tiny
\begin{split}
&\;\;\;\; ( \overline{\partial_1}+\mathcal{D}_1)\left(e^{i_1}e^{i_2}\cdots e^{i_r}\right)\\
&=\left[ \delta\binom{i_1}{i_2\epsilon}+2\delta\binom{i_1}{-i_2}+\delta\binom{i_1}{i_2\epsilon^{-1}}+\delta\binom{i_1}{i_2}\right]e^{i_2}\cdots e^{i_r}+\cdots+\left[\delta\binom{i_{r-1}}{i_r\epsilon}+2\delta\binom{i_{r-1}}{-i_r}+\delta\binom{i_{r-1}}{i_r\epsilon^{-1}} +\delta\binom{i_{r-1}}{i_r}\right]e^{i_1}\cdots e^{i_{r-2}}e^{i_r}    \\
&+\left[\delta\binom{i_r}{\epsilon}+2\delta\binom{i_r}{-1}+\delta\binom{i_r}{\epsilon^{-1}}+\delta\binom{i_r}{1}  \right]e^{i_1}\cdots e^{i_{r-1}}-\left[\delta\binom{i_1}{i_2\epsilon}+2\delta\binom{i_1}{-i_2}+\delta\binom{i_1}{i_2\epsilon^{-1}}+\delta\binom{i_1}{i_2} \right]e^{i_1}e^{i_3}\cdots e^{i_r}\\
&-\cdots- \left[\delta\binom{i_{r-1}}{i_r\epsilon}+2\delta\binom{i_{r-1}}{-i_r}+\delta\binom{i_{r-1}}{i_r\epsilon^{-1}} +\delta\binom{i_{r-1}}{i_r}\right]e^{i_1}\cdots e^{i_{r-2}}e^{i_{r-1}}.  \\
&=\left[ 1+\delta\binom{i_1}{-i_2}\right]e^{i_2}\cdots e^{i_r}+\cdots+\left[1+\delta\binom{i_{r-1}}{-i_r}\right]e^{i_1}\cdots e^{i_{r-2}}e^{i_r}  +\left[1+\delta\binom{i_r}{-1}\right]e^{i_1}\cdots e^{i_{r-1}}\\
&\;\;\;\;-\left[1+\delta\binom{i_1}{-i_2} \right]e^{i_1}e^{i_3}\cdots e^{i_r}
-\cdots- \left[1+\delta\binom{i_{r-1}}{-i_r}\right]e^{i_1}\cdots e^{i_{r-2}}e^{i_{r-1}} \\
&=e^{i_2}\cdots e^{i_r}+\delta\binom{i_1}{-i_2}e^{i_2}\cdots e^{i_r}+\cdots+\delta\binom{i_{r-1}}{-i_r}e^{i_1}\cdots e^{i_{r-2}}e^{i_r}  +\delta\binom{i_r}{-1}e^{i_1}\cdots e^{i_{r-1}}\\
&\;\;\;\;-\delta\binom{i_1}{-i_2} e^{i_1}e^{i_3}\cdots e^{i_r}
-\cdots- \delta\binom{i_{r-1}}{-i_r}e^{i_1}\cdots e^{i_{r-2}}e^{i_{r-1}} \\
&=(\mathcal{D}_2+\mathcal{P})\left(e^{i_1}e^{i_2}\cdots e^{i_r}\right).
\end{split}
\]
$(vi)$ This statement follows from the following observation
\[
\delta\binom{\eta}{-i_1}=0, \;\mathrm{if}\; \eta/i_1\in \{\epsilon^{\pm 1}\}.
\]
$\hfill\Box$\\

Now we are ready to prove Theorem \ref{efca}.

For $(i)$ and $(ii)$, $N=2,3 $, by Theorem \ref{mmm} and Lemma \ref{below}, it suffices to show that
\[
\underbrace{\overline{\partial_1}\circ\cdots \circ \overline{\partial_1}}_{r-1}\left(  e^{\eta_1}\cdots e^{\eta_r}  \right)=e^{\eta_r},\,\eta_1,\cdots,\eta_r\in \mu_N,\,\eta_r\neq 1, \,r\geq 2.
\]
By Proposition \ref{pd}, $(i)$ and $(iv)$,  for $\eta_r\neq 1$, one has
\[
\overline{\partial_1}(e^{\eta_1}\cdots e^{\eta_r})=(\mathcal{D}_2-\mathcal{D}_1)(e^{\eta_1}\cdots e^{\eta_r})=e^{\eta_2}\cdots e^{\eta_r}.
\]
By induction, it follows that
\[
\underbrace{\overline{\partial_1}\circ\cdots \circ \overline{\partial_1}}_{r-1}\left(  e^{\eta_1}\cdots e^{\eta_r}  \right)=e^{\eta_r}.
\]
As a result, the statements $(i)$ and $(ii)$ are proved.

For $(iii)$, $N=4$, by Lemma \ref{below}, one has
\[
\underbrace{\overline{\partial_1}\circ \cdots \circ\overline{\partial_1}}_{r-1}(\underbrace{e^{-1}\cdots e^{-1}}_r )=2^{r-1}e^{-1},
\]
\[
\underbrace{\overline{\partial_1}\circ \cdots \circ\overline{\partial_1}}_{r-1}(\underbrace{e^{\eta}\cdots e^{\eta}}_r )=e^{\eta}, \eta\in \{\epsilon, \epsilon^{-1}\},
\]
\[
\overline{\partial_1}(e^{-1}-e^{\epsilon}-e^{\epsilon^{-1}})=0.
\]
Thus by Theorem \ref{mmm}, we have
\[
\zeta\binom{1,\cdots,1,\;\;1}{\underbrace{1,\cdots,1,-1}_r }=2^{r-1}\left(   \zeta\binom{1,\cdots,1,1}{\underbrace{1,\cdots,1,\epsilon}_r }+    \zeta\binom{1,\cdots,1,\;1}{\underbrace{1,\cdots,1,\epsilon^{-1}}_r }     \right)+\mathrm{LDT}.\tag{1}\]

For the formula $(2)$, it suffices to show that
\[
\underbrace{\overline{\partial_1}\circ \cdots \circ\overline{\partial_1}}_{s-1} (\underbrace{e^{\eta_1}\cdots e^{\eta_{s-1}}e^{\eta_s} } )=2^{s-1}e^{\eta_s}, \;\eta_1,\cdots, \eta_{r-1}\in \mu_2,\;\eta_s=-1.
\]
By definition, for $N=4$, $\eta_1,\cdots, \eta_{s-1}\in \mu_2,\;\eta_s=-1$, one has
\[
\overline{\partial_1}(e^{\eta_1}\cdots e^{\eta_{s-1}}e^{\eta_s})=2\mathcal{P}(e^{\eta_1}\cdots e^{\eta_{s-1}}e^{\eta_s}).\]
By the same analysis in the proof of $(i)$ and $(ii)$, it follows that 
\[
\underbrace{\overline{\partial_1}\circ\cdots \circ \overline{\partial_1}}_{r-1}(e^{\eta_1}\cdots e^{\eta_{r-1}}e^{\eta_r})=2^{r-1}e^{\eta_r},\;\eta_1,\cdots, \eta_{r-1}\in \mu_2,\;\eta_r=-1.\]
Thus 
\[\zeta\binom{1,\;1,\;\cdots,1}{\eta_1,\eta_2,\cdots, \eta_s}=\zeta\binom{1,\cdots,1,\,1}{1,\cdots,1,\eta_s}+\mathrm{LDT},   \eta_j \in \{\pm 1\}, \forall\, j ,\eta_s=-1.\tag{2}
\]

For the formula $(3)$, it suffices to show that 
\[
\underbrace{\overline{\partial_1}\circ \cdots \circ \overline{\partial_1}}_{r-1}(e^{\epsilon_1\cdots \epsilon_r}\cdots e^{\epsilon_{r-1}\epsilon_r}e^{\epsilon_r})=e^{\epsilon_r}, \epsilon_1,\cdots, \epsilon_r\in \{\epsilon, \epsilon^{-1}\}.
\]
By Proposition \ref{pd}, $(v)$, for $\epsilon_1,\cdots, \epsilon_r\in \{\epsilon, \epsilon^{-1}\}$, one has
\[
\begin{split}
&\;\;\;\;\,\overline{\partial_1}\left(e^{\epsilon_1\cdots \epsilon_r}\cdots e^{\epsilon_{r-1}\epsilon_r}e^{\epsilon_r}\right)\\
&=\left(\mathcal{D}_2+\mathcal{P}-\mathcal{D}_1\right)\left(  e^{\epsilon_1\cdots \epsilon_r}\cdots e^{\epsilon_{r-1}\epsilon_r}e^{\epsilon_r}   \right)\\
&=\mathcal{D}_2\left(  e^{\epsilon_1\cdots \epsilon_r}\cdots e^{\epsilon_{r-1}\epsilon_r}e^{\epsilon_r}   \right)\\
&= e^{\epsilon_2\cdots \epsilon_r}\cdots e^{\epsilon_{r-1}\epsilon_r}e^{\epsilon_r} .
\end{split}\]
By induction, it follows that
\[
\underbrace{\overline{\partial_1}\circ \cdots \circ \overline{\partial_1}}_{r-1}(e^{\epsilon_1\cdots \epsilon_r}\cdots e^{\epsilon_{r-1}\epsilon_r}e^{\epsilon_r})=e^{\epsilon_r}, \epsilon_1,\cdots, \epsilon_r\in \{\epsilon, \epsilon^{-1}\}.\]
Thus 
\[
\zeta\binom{1,\;1,\;\cdots,1}{\epsilon_1,\epsilon_2,\cdots, \epsilon_r}=
\zeta\binom{1,\cdots,1,\,1}{1,\cdots,1,\epsilon_r}+\mathrm{LDT}, \epsilon_j \in \{\epsilon,\epsilon^{-1}\}, \forall\,j.\tag{3}
\]

For the formula $(4)$, it suffices to show that: For $\epsilon_1,\cdots, \epsilon_r\in \{\epsilon,\epsilon^{-1}\}, \eta_1,\cdots, \eta_s\in \{\pm 1\}, \eta_1=-1, r,s\geq 1$, one has
\[
\begin{split}
&\;\;\;\; \underbrace{\overline{\partial_1}\circ \cdots \circ \overline{\partial_1}}_{r+s-1}\left(  e^{\epsilon_1\cdots \epsilon_r\eta_1\cdots \eta_s}\cdots e^{\epsilon_1\eta_1\cdots \eta_s} e^{\eta_1\cdots \eta_s}  \cdots e^{\eta_1} \right)\\
&= \sum_{j=0}^{s-1}\binom{r+s-1}{j}\cdot e^{\eta_1} +\binom{r+s-1}{s} e^{\epsilon_1\eta_1\cdots \eta_s}.    \\
\end{split}
\]
By Proposition \ref{pd}, $(v)$, one has
\[
\begin{split}
&\;\;\;\; \overline{\partial_1}\left(  e^{\epsilon_1\cdots \epsilon_r\eta_1\cdots \eta_s}\cdots e^{\epsilon_1\eta_1\cdots \eta_s} e^{\eta_1\cdots \eta_s}  \cdots e^{\eta_1} \right)\\
&=\left( \overline{\partial_1} +\mathcal{D}_1  \right)\left(  e^{\epsilon_1\cdots \epsilon_r\eta_1\cdots \eta_s}\cdots e^{\epsilon_1\eta_1\cdots \eta_s} e^{\eta_1\cdots \eta_s}  \cdots e^{\eta_1} \right) \\
&=\left(\mathcal{D}_2+\mathcal{P}\right)\left(  e^{\epsilon_1\cdots \epsilon_r\eta_1\cdots \eta_s}\cdots e^{\epsilon_1\eta_1\cdots \eta_s} e^{\eta_1\cdots \eta_s}  \cdots e^{\eta_1} \right) .\\
\end{split}
\]
By Proposition \ref{pd}, $(iii)$, we have
\[
\begin{split}
&\;\;\;\;  \underbrace{(\overline{\partial_1}+\mathcal{D}_1)\circ \cdots \circ (\overline{\partial_1}+\mathcal{D}_1) }_{r+s-1} \left(  e^{\epsilon_1\cdots \epsilon_r\eta_1\cdots \eta_s}\cdots e^{\epsilon_1\eta_1\cdots \eta_s} e^{\eta_1\cdots \eta_s}  \cdots e^{\eta_1} \right)      \\
&=\left[ \left(  \overline{\partial_1} \right)^{r+s-1}  +\sum_{j=1}^{r+s-1} \binom{r+s-1}{j}  \left(  \overline{\partial_1} \right)^{r+s-j-1}  (\mathcal{D}_1)^j \right]  \left(  e^{\epsilon_1\cdots \epsilon_r\eta_1\cdots \eta_s}\cdots e^{\epsilon_1\eta_1\cdots \eta_s} e^{\eta_1\cdots \eta_s}  \cdots e^{\eta_1} \right)      \\
&= \underbrace{\overline{\partial_1}\circ \cdots \circ \overline{\partial_1} }_{r+s-1} \left(  e^{\epsilon_1\cdots \epsilon_r\eta_1\cdots \eta_s}\cdots e^{\epsilon_1\eta_1\cdots \eta_s} e^{\eta_1\cdots \eta_s}  \cdots e^{\eta_1} \right)      \\
&=\sum_{j=0}^{r+s-1}\binom{r+s-1}{j}  \mathcal{P}^j\circ \mathcal{D}_2^{r+s-1-j}  \left(  e^{\epsilon_1\cdots \epsilon_r\eta_1\cdots \eta_s}\cdots e^{\epsilon_1\eta_1\cdots \eta_s} e^{\eta_1\cdots \eta_s}  \cdots e^{\eta_1} \right)             \\
&=\sum_{j=0}^{s-1}\binom{r+s-1}{j}  \mathcal{P}^j\left( e^{\eta_1\cdots \eta_{j+1}}  \cdots e^{\eta_1} \right)\\
&\;\;\;\;+\sum_{j=s}^{r+s-1}\binom{r+s-1}{j}  \mathcal{P}^j\left(  e^{\epsilon_1\cdots \epsilon_{j-s+1}\eta_1\cdots \eta_s} \cdots  \cdots e^{\epsilon_1\eta_1\cdots \eta_s}e^{\eta_1\cdots \eta_{s}}  \cdots e^{\eta_1} \right).
\end{split}
\]
By Proposition \ref{pd}, $(vi)$, it follows that:\\
For $j\geq s+1$,
\[
\begin{split}
&\;\;\;\;\mathcal{P}^j\left(  e^{\epsilon_1\cdots \epsilon_{j-s+1}\eta_1\cdots \eta_s} \cdots  \cdots e^{\epsilon_1\eta_1\cdots \eta_s}e^{\eta_1\cdots \eta_{s}}  \cdots e^{\eta_1} \right)\\
&=e^{\epsilon_1\cdots \epsilon_{j-s+1}\eta_1\cdots \eta_s} \mathcal{P}^j\left(  e^{\epsilon_1\cdots \epsilon_{j-s}\eta_1\cdots \eta_s} \cdots  \cdots e^{\epsilon_1\eta_1\cdots \eta_s}e^{\eta_1\cdots \eta_{s}}  \cdots e^{\eta_1} \right)  \\
&\cdots\\
&=0.
\end{split}
\]
For $j=s$, 
\[
\begin{split}
&\;\;\;\;\mathcal{P}^s\left(   e^{\epsilon_1\eta_1\cdots \eta_s}e^{\eta_1\cdots \eta_{s}}  \cdots e^{\eta_1} \right)\\
&=e^{\epsilon_1\eta_1\cdots \eta_s} \mathcal{P}^s\left(  e^{\eta_1\cdots \eta_{s}}  \cdots e^{\eta_1} \right).  \\
\end{split}
\]
By the same analysis in the proof of $(i)$ and $(ii)$, one has 
\[
 \mathcal{P}^j\left( e^{\eta_1\cdots \eta_{j+1}}  \cdots e^{\eta_1} \right)=e^{\eta_1}, j=0,\cdots, {s-1},
 \]
 \[
  \mathcal{P}^s\left(  e^{\eta_1\cdots \eta_{s}}  \cdots e^{\eta_1} \right)=1. \]
  As a result, \[
\begin{split}
&\;\;\;\; \underbrace{\overline{\partial_1}\circ \cdots \circ \overline{\partial_1}}_{r+s-1}\left(  e^{\epsilon_1\cdots \epsilon_r\eta_1\cdots \eta_s}\cdots e^{\epsilon_1\eta_1\cdots \eta_s} e^{\eta_1\cdots \eta_s}  \cdots e^{\eta_1} \right)\\
&= \sum_{j=0}^{s-1}\binom{r+s-1}{j}\cdot e^{\eta_1} +\binom{r+s-1}{s} e^{\epsilon_1\eta_1\cdots \eta_s}.    \\
\end{split}
\]

For the formula $(5)$, it suffices to show that : \\
For $s\geq 1$, $\epsilon_1,\cdots, \epsilon_{s+1}\in \{\epsilon,\epsilon^{-1}=-\epsilon\}$, one has 
\[
\underbrace{\overline{\partial_1}\circ \cdots \circ \overline{\partial_1}}_{s} \left(e^{\epsilon_{s+1}}e^{\epsilon_s}\cdots e^{\epsilon_1} \right)=\sum_{j=0}^s(-1)^j\binom{s}{j} 2^{s-j}e^{\epsilon_{j+1}}.
\]
 For $\epsilon_1,\cdots, \epsilon_{s+1}\in \{\epsilon,\epsilon^{-1}=-\epsilon\}$,  we have
 \[
 \begin{split}
 &\;\;\;\;   \mathcal{P}\left(  e^{\epsilon_{s+1}}e^{\epsilon_s}\cdots e^{\epsilon_1}    \right)           \\
 &=\delta\binom{\epsilon_{s+1}}{-\epsilon_s} e^{\epsilon_s} e^{\epsilon_{s-1}}\cdots e^{\epsilon_1} +\cdots +\delta\binom{\epsilon_2}{-\epsilon_1}   e^{\epsilon_{s+1}}\cdots e^{\epsilon_3} e^{\epsilon_1}        \\
  &\;\;-\delta\binom{\epsilon_{s+1}}{-\epsilon_s} e^{\epsilon_{s+1}} e^{\epsilon_{s-1}} \cdots e^{\epsilon_1} +\cdots +\delta\binom{\epsilon_2}{-\epsilon_1}   e^{\epsilon_{s+1}}\cdots e^{\epsilon_3} e^{\epsilon_2}        \\
  &=\left[\delta\binom{\epsilon_{s+1}}{-\epsilon_s}+\delta\binom{\epsilon_{s+1}}{\epsilon_s} \right] e^{\epsilon_s} e^{\epsilon_{s-1}}\cdots e^{\epsilon_1} +\cdots +\left[\delta\binom{\epsilon_2}{-\epsilon_1}+\delta\binom{\epsilon_2}{\epsilon_1} \right]  e^{\epsilon_{s+1}}\cdots e^{\epsilon_3} e^{\epsilon_1}        \\
  &\;\;-\left[\delta\binom{\epsilon_{s+1}}{-\epsilon_s}+\delta\binom{\epsilon_{s+1}}{\epsilon_s} \right] e^{\epsilon_{s+1}} e^{\epsilon_{s-1}} \cdots e^{\epsilon_1} -\cdots -\left[\delta\binom{\epsilon_2}{-\epsilon_1}+\delta\binom{\epsilon_2}{\epsilon_1} \right]  e^{\epsilon_{s+1}}\cdots e^{\epsilon_3} e^{\epsilon_2}        \\ 
  &=e^{\epsilon_s} e^{\epsilon_{s-1}}\cdots e^{\epsilon_1} + e^{\epsilon_{s+1}} e^{\epsilon_{s-1}} \cdots e^{\epsilon_1}+\cdots +e^{\epsilon_{s+1}}\cdots e^{\epsilon_3} e^{\epsilon_1}        \\
  &\;\;- e^{\epsilon_{s+1}} e^{\epsilon_{s-1}} \cdots e^{\epsilon_1} -\cdots -e^{\epsilon_{s+1}}\cdots e^{\epsilon_3} e^{\epsilon_1}  - e^{\epsilon_{s+1}}\cdots e^{\epsilon_3} e^{\epsilon_2}        \\ 
  &= e^{\epsilon_s} e^{\epsilon_{s-1}}\cdots e^{\epsilon_1}-  e^{\epsilon_{s+1}}\cdots e^{\epsilon_3} e^{\epsilon_2}\\
  &=\mathcal{D}_2\left(  e^{\epsilon_{s+1}}e^{\epsilon_s}\cdots e^{\epsilon_1}    \right)-\widetilde{\mathcal{D}}_1 \left(  e^{\epsilon_{s+1}}e^{\epsilon_s}\cdots e^{\epsilon_1}    \right),   \end{split}
 \]
 where the linear map $\widetilde{\mathcal{D}}_1: \mathbb{Q}\langle e^{\mu_N}\rangle \rightarrow  \mathbb{Q}\langle e^{\mu_N}\rangle$ is defined by 
 \[
 \widetilde{\mathcal{D}}_1(1)=0, \; \widetilde{\mathcal{D}}_1(e^{i_1})=1,\;\widetilde{\mathcal{D}}_1(e^{i_r}\cdots e^{i_2} e^{i_1})=e^{i_r}\cdots e^{i_2} .\]
  Thus 
 \[
 \begin{split}
 &\;\;\;\; \,\overline{\partial_1} \left(  e^{\epsilon_{s+1}}e^{\epsilon_s}\cdots e^{\epsilon_1}    \right)        \\
  &= \left(\overline{\partial_1} +\mathcal{D}_1       \right)  \left( e^{\epsilon_{s+1}}e^{\epsilon_s}\cdots e^{\epsilon_1}    \right)            \\
 &= \left(\mathcal{D}_2+\mathcal{P}    \right)  \left( e^{\epsilon_{s+1}}e^{\epsilon_s}\cdots e^{\epsilon_1}    \right)            \\
 &=\left(2\mathcal{D}_2-  \widetilde{\mathcal{D}}_1   \right)\left( e^{\epsilon_{s+1}}e^{\epsilon_s}\cdots e^{\epsilon_1}    \right).             \\
 \end{split}
 \]
 Since it is clear that $\mathcal{D}_2 \widetilde{\mathcal{D}}_1 =\widetilde{\mathcal{D}}_1 \mathcal{D}_2$, we have 
 \[
 \begin{split}
&\;\;\;\;\underbrace{\overline{\partial_1}\circ \cdots \circ \overline{\partial_1}}_{s} \left(e^{\epsilon_{s+1}}e^{\epsilon_s}\cdots e^{\epsilon_1} \right)\\
&=\sum_{j=0}^s(-1)^j\binom{s}{j} 2^{s-j} \left(\mathcal{D}_2\right)^{s-j}  \circ  \left(\widetilde{\mathcal{D}}_1\right)^j \left(e^{\epsilon_{s+1}}e^{\epsilon_s}\cdots e^{\epsilon_1} \right)     \\
&=\sum_{j=0}^s(-1)^j\binom{s}{j} 2^{s-j} e^{\epsilon_{j+1}}.
\end{split}
 \]
 
 For the formula $(6)$, it suffices to show that:\\
 For $\epsilon_1,\cdots , \epsilon_r\in \{\epsilon,\epsilon^{-1}\}, \eta_1,\cdots, \eta_s\in \{\pm 1\}, r\geq 2, s\geq 1$, 
 \[
\begin{split}
&\;\;\;\; \underbrace{\overline{\partial_1}\circ \cdots \circ \overline{\partial_1}}_{r+s-1}\left(  e^{\eta_1\cdots \eta_s\epsilon_1\cdots \epsilon_r}\cdots e^{\eta_1\epsilon_1\cdots \epsilon_r} e^{\epsilon_1\cdots \epsilon_r}  \cdots e^{\epsilon_1} \right)=e^{\epsilon_1}.
\end{split} \]
By Proposition \ref{pd}, $(iii)$, we have
\[
\begin{split}
&\;\;\;\; \underbrace{\overline{\partial_1}\circ \cdots \circ \overline{\partial_1}}_{r+s-1}\left(  e^{\eta_1\cdots \eta_s\epsilon_1\cdots \epsilon_r}\cdots e^{\eta_1\epsilon_1\cdots \epsilon_r} e^{\epsilon_1\cdots \epsilon_r}  \cdots e^{\epsilon_1} \right)\\
&= \underbrace{(\overline{\partial_1}+\mathcal{D}_1)\circ \cdots \circ (\overline{\partial_1}+\mathcal{D}_1) }_{r+s-1} \left(   e^{\eta_1\cdots \eta_s\epsilon_1\cdots \epsilon_r}\cdots e^{\eta_1\epsilon_1\cdots \epsilon_r} e^{\epsilon_1\cdots \epsilon_r}  \cdots e^{\epsilon_1} \right)      \\
&=  \underbrace{(\mathcal{D}_2+\mathcal{P})\circ \cdots \circ (\mathcal{D}_2+\mathcal{P}) }_{r+s-1} \left(  e^{\eta_1\cdots \eta_s\epsilon_1\cdots \epsilon_r}\cdots e^{\eta_1\epsilon_1\cdots \epsilon_r} e^{\epsilon_1\cdots \epsilon_r}  \cdots e^{\epsilon_1} \right)       \\
&=\left[ \left(\mathcal{D}_2\right)^{r+s-1}+\sum_{j=1}^{r+s-1} \binom{r+s-1}{j} \mathcal{P}^{j}\circ \left(\mathcal{D}_2\right)^{r+s-1-j} \right] \left(  e^{\eta_1\cdots \eta_s\epsilon_1\cdots \epsilon_r}\cdots e^{\eta_1\epsilon_1\cdots \epsilon_r} e^{\epsilon_1\cdots \epsilon_r}  \cdots e^{\epsilon_1} \right)     \\
&=e^{\epsilon_1}+\sum_{j=1}^{r-1} \binom{r+s-1}{j} \mathcal{P}^{j}\left(e^{\epsilon_1\cdots \epsilon_{j+1}}\cdots e^{\epsilon_1}   \right)\\
&\;\;\;\;+\sum_{j=r}^{r+s-1} \binom{r+s-1}{j} \mathcal{P}^{j}
 \left(  e^{\eta_1\cdots \eta_{j-r+1}\epsilon_1\cdots \epsilon_r}\cdots e^{\eta_1\epsilon_1\cdots \epsilon_r} e^{\epsilon_1\cdots \epsilon_r}  \cdots e^{\epsilon_1} \right)\\  
 &=e^{\epsilon_1}+\sum_{j=r}^{r+s-1} \binom{r+s-1}{j} \mathcal{P}^{j}
 \left(  e^{\eta_1\cdots \eta_{j-r+1}\epsilon_1\cdots \epsilon_r}\cdots e^{\eta_1\epsilon_1\cdots \epsilon_r} e^{\epsilon_1\cdots \epsilon_r}  \cdots e^{\epsilon_1} \right).   
  \end{split}
\] 
As a result, the formula $(6)$ is reduced to the following statement:\\
  For $$\epsilon_1,\cdots , \epsilon_r\in \{\epsilon,\epsilon^{-1}\}, \eta_1,\cdots, \eta_s\in \{\pm 1\},$$ $$ r\geq 2, s\geq 1, r\leq j\leq r+s-1,$$
   one has 
   \[
   \mathcal{P}^{j}
 \left(  e^{\eta_1\cdots \eta_{j-r+1}\epsilon_1\cdots \epsilon_r}\cdots e^{\eta_1\epsilon_1\cdots \epsilon_r} e^{\epsilon_1\cdots \epsilon_r}  \cdots e^{\epsilon_1} \right)=0.  
   \]
   We have
   \[\tiny
   \begin{split}
   &\;\;\;\;  \mathcal{P}    \left(  e^{\eta_1\cdots \eta_{j-r+1}\epsilon_1\cdots \epsilon_r}\cdots e^{\eta_1\epsilon_1\cdots \epsilon_r} e^{\epsilon_1\cdots \epsilon_r}  \cdots e^{\epsilon_1} \right)                 \\
   &=  \delta\binom{\eta_{j-r+1}}{-1}  e^{\eta_1\cdots \eta_{j-r}\epsilon_1\cdots \epsilon_r}\cdots e^{\eta_1\epsilon_1\cdots \epsilon_r} e^{\epsilon_1\cdots \epsilon_r}  \cdots e^{\epsilon_1} +\cdots +\delta\binom{\eta_1}{-1}   e^{\eta_1\cdots \eta_{j-r+1}\epsilon_1\cdots \epsilon_r}\cdots e^{\eta_1\eta_2\epsilon_1\cdots \epsilon_r} e^{\epsilon_1\cdots \epsilon_r}  \cdots e^{\epsilon_1}                        \\
   &\;\;\;\;-     \delta\binom{\eta_{j-r+1}}{-1}  e^{\eta_1\cdots \eta_{j-r+1}\epsilon_1\cdots \epsilon_r} e^{\eta_1\cdots \eta_{j-r-1}\epsilon_1\cdots \epsilon_r}   \cdots e^{\eta_1\epsilon_1\cdots \epsilon_r} e^{\epsilon_1\cdots \epsilon_r}  \cdots e^{\epsilon_1} -\cdots \\
   &\;\;\;\;-\delta\binom{\eta_1}{-1}   e^{\eta_1\cdots \eta_{j-r+1}\epsilon_1\cdots \epsilon_r}\cdots e^{\eta_1\epsilon_1\cdots \epsilon_r} e^{\epsilon_1\cdots \epsilon_{r-1}}  \cdots e^{\epsilon_1}                                   \\
   &=\left[  \delta\binom{\eta_{j-r+1}}{-1} + \delta\binom{\eta_{j-r+1}}{1}  \right] e^{\eta_1\cdots \eta_{j-r}\epsilon_1\cdots \epsilon_r}\cdots e^{\eta_1\epsilon_1\cdots \epsilon_r} e^{\epsilon_1\cdots \epsilon_r}  \cdots e^{\epsilon_1}+\cdots   \\
   &\;\;\;\;+ \left[  \delta\binom{\eta_{1}}{-1} + \delta\binom{\eta_{1}}{1}  \right] e^{\eta_1\cdots \eta_{j-r+1}\epsilon_1\cdots \epsilon_r}\cdots e^{\eta_1\eta_2\epsilon_1\cdots \epsilon_r} e^{\epsilon_1\cdots \epsilon_r}  \cdots e^{\epsilon_1}  \\
    &\;\;\;\;-     \left[ \delta\binom{\eta_{j-r+1}}{-1} + \delta\binom{\eta_{j-r+1}}{1}    \right]  e^{\eta_1\cdots \eta_{j-r+1}\epsilon_1\cdots \epsilon_r} e^{\eta_1\cdots \eta_{j-r-1}\epsilon_1\cdots \epsilon_r}   \cdots e^{\eta_1\epsilon_1\cdots \epsilon_r} e^{\epsilon_1\cdots \epsilon_r}  \cdots e^{\epsilon_1} -\cdots \\
    &\;\;\;\;-\left[ \delta\binom{\eta_1}{-1}  + \delta\binom{\eta_1}{1}  \right]   e^{\eta_1\cdots \eta_{j-r+1}\epsilon_1\cdots \epsilon_r}\cdots e^{\eta_1\epsilon_1\cdots \epsilon_r} e^{\epsilon_1\cdots \epsilon_{r-1}}  \cdots e^{\epsilon_1}                                   \\   
    &= e^{\eta_1\cdots \eta_{j-r}\epsilon_1\cdots \epsilon_r}\cdots e^{\eta_1\epsilon_1\cdots \epsilon_r} e^{\epsilon_1\cdots \epsilon_r}  \cdots e^{\epsilon_1}+\cdots +   e^{\eta_1\cdots \eta_{j-r+1}\epsilon_1\cdots \epsilon_r}\cdots e^{\eta_1\eta_2\epsilon_1\cdots \epsilon_r} e^{\epsilon_1\cdots \epsilon_r}  \cdots e^{\epsilon_1}\\
    &\;\;\;\;- e^{\eta_1\cdots \eta_{j-r+1}\epsilon_1\cdots \epsilon_r} e^{\eta_1\cdots \eta_{j-r-1}\epsilon_1\cdots \epsilon_r}   \cdots e^{\eta_1\epsilon_1\cdots \epsilon_r} e^{\epsilon_1\cdots \epsilon_r}  \cdots e^{\epsilon_1} -\cdots -      e^{\eta_1\cdots \eta_{j-r+1}\epsilon_1\cdots \epsilon_r}\cdots e^{\eta_1\epsilon_1\cdots \epsilon_r} e^{\epsilon_1\cdots \epsilon_{r-1}}  \cdots e^{\epsilon_1}                            \\
    &= e^{\eta_1\cdots \eta_{j-r}\epsilon_1\cdots \epsilon_r}\cdots e^{\eta_1\epsilon_1\cdots \epsilon_r} e^{\epsilon_1\cdots \epsilon_r}  \cdots e^{\epsilon_1}-    e^{\eta_1\cdots \eta_{j-r+1}\epsilon_1\cdots \epsilon_r}\cdots e^{\eta_1\epsilon_1\cdots \epsilon_r} e^{\epsilon_1\cdots \epsilon_{r-1}}  \cdots e^{\epsilon_1}  .    \\
                          \end{split}
   \]
   Since $\eta_1\epsilon_r \in \{\epsilon, \epsilon^{-1}\}$, by induction, 
   \[
   \mathcal{P}^{j-r+1}    \left(  e^{\eta_1\cdots \eta_{j-r+1}\epsilon_1\cdots \epsilon_r}\cdots e^{\eta_1\epsilon_1\cdots \epsilon_r} e^{\epsilon_1\cdots \epsilon_r}  \cdots e^{\epsilon_1} \right)     \]
   is a $\mathbb{Q}$-linear combination of the following kinds of elements:
   \[
   e^{i_1\cdots i_r}\cdots e^{i_1}, i_1,\cdots,i_r\in \{\epsilon,\epsilon^{-1}\}.
   \]
   For $r\geq 2$, by the definition of $\mathcal{P}$, it follows that
   \[
   \begin{split}
   & \;\;\;\;\,\mathcal{P}^{j}
 \left(  e^{\eta_1\cdots \eta_{j-r+1}\epsilon_1\cdots \epsilon_r}\cdots e^{\eta_1\epsilon_1\cdots \epsilon_r} e^{\epsilon_1\cdots \epsilon_r}  \cdots e^{\epsilon_1} \right)\\
 &= \mathcal{P}^{r-1} \circ  \mathcal{P}^{j-r+1}
 \left(  e^{\eta_1\cdots \eta_{j-r+1}\epsilon_1\cdots \epsilon_r}\cdots e^{\eta_1\epsilon_1\cdots \epsilon_r} e^{\epsilon_1\cdots \epsilon_r}  \cdots e^{\epsilon_1} \right) \\
 &=0.
 \end{split}   \]
 Thus the formula $(6)$ is proved.
 
 \begin{ex}
 For $N=4$, $\mu_4=\{1,-1,i,-i\}$. One has 
 \[
  \zeta\binom{1}{-1}=\sum_{n\geq 1}\frac{(-1)^n}{n}=-\mathrm{log}\,2, \]
 \[
 \zeta\binom{1}{i}=\sum_{n\geq 1}\frac{i^n}{n}=-\mathrm{log}\,(1-i)=-\frac{\mathrm{log}\,2}{2}+\frac{\pi}{4}i,
 \]
 \[
 \zeta\binom{1}{-i}=\sum_{n\geq 1}\frac{(-i)^n}{n}=-\mathrm{log}\,(1+i)=-\frac{\mathrm{log}\,2}{2}-\frac{\pi}{4}i.
 \]
 For $r\geq 2$, it follows that
 \[
 \begin{split}
 &\;\;\;\;   \left[\zeta\binom{1}{i} \right]^r+ \left[ \zeta\binom{1}{-i}\right]^r\\
 &= \left(    -\frac{\mathrm{log}\,2}{2}+\frac{\pi}{4}i   \right)^r+  \left(   -\frac{\mathrm{log}\,2}{2}-\frac{\pi}{4}i    \right)^r          \\
 &=2\sum_{\substack{0\leq j\leq r\\ j\,\equiv \,0\,\mathrm{mod}\, 2}}\binom{r}{j}\left( -\frac{\mathrm{log}\,2}{2}  \right)^{r-j}\left(\frac{\pi }{4} i\right)^j\\
 &= \frac{1}{2^{r-1}} \left(-\mathrm{log}\,2\right)^r +\frac{1}{2^{r+j-1}} \sum_{\substack{0\leq j\leq r\\ j\,\equiv \,0\,\mathrm{mod}\, 2}}  (-1)^{\frac{j}{2}}\binom{r}{j}\left(-\mathrm{log}\,2\right)^{r-j} \pi^j.
  \end{split}
 \]
 Thus 
 \[
 \begin{split}
 &\;\;\;\;   \zeta\binom{1,\cdots,1,\;\,1\;}{\underbrace{1,\cdots,1,-1}_r }     \\
 &= \frac{1}{r!}  \left(  -\mathrm{log}\,2   \right)^r           \\
  &=\frac{2^{r-1}}{r!}    \left[\zeta\binom{1}{i} \right]^r+ \frac{2^{r-1}}{r!}  \left[ \zeta\binom{1}{-i}\right]^r - \frac{1}{2^{j}} \sum_{\substack{0\leq j\leq r\\ j\,\equiv \,0\,\mathrm{mod}\, 2}}  (-1)^{\frac{j}{2}}\frac{1}{j!(r-j)!}\left(-\mathrm{log}\,2\right)^{r-j} \pi^j   \\
  &=2^{r-1}\left(   \zeta\binom{1,\cdots,1,1}{\underbrace{1,\cdots,1,\epsilon}_r }+    \zeta\binom{1,\cdots,1,\;1}{\underbrace{1,\cdots,1,\epsilon^{-1}}_r }     \right) -\frac{1}{2^{j}} \sum_{\substack{0\leq j\leq r\\ j\,\equiv \,0\,\mathrm{mod}\, 2}}  (-1)^{\frac{j}{2}}   \zeta\binom{1,\cdots,1,\;\,1\;}{\underbrace{1,\cdots,1,-1}_{r-j} }\cdot \frac{\pi^j }{j!} .     \end{split}
 \]
  Since $\zeta(j)/\pi^j\in\mathbb{Q}$ for $r\geq 2$, even, by the stuffle product of cyclotomic multiple zeta values,  the element \[\sum_{\substack{0\leq j\leq r\\ j\,\equiv \,0\,\mathrm{mod}\, 2}}  (-1)^{\frac{j}{2}}   \zeta\binom{1,\cdots,1,\;\,1\;}{\underbrace{1,\cdots,1,-1}_{r-j} }\cdot \frac{\pi^j }{j!}  \] belongs to the $\mathbb{Q}$-linear combinations of cyclotomic multiple zeta values of depth $<r$. This is compatible with the formula $(1)$ in Theorem \ref{gene}.
   \end{ex}
 
 \section*{Acknowledgements}
The author is supported by the National Natural Science Foundation of China (Grant No. 12201642).

\end{document}